\documentclass{amsart}

\usepackage{enumerate}

\numberwithin{equation}{section}

\newcommand{\bF}{\mathbb{F}}
\newcommand{\bG}{\mathbb{G}}

\newcommand{\bN}{\mathbb{N}}
\newcommand{\bP}{\mathbb{P}}
\newcommand{\bQ}{\mathbb{Q}}
\newcommand{\bR}{\mathbb{R}}
\newcommand{\bS}{\mathbb{S}}
\newcommand{\bT}{\mathbb{T}}
\newcommand{\bZ}{\mathbb{Z}}

\newcommand{\bone}{\mathbf{1}}

\newcommand{\cD}{\mathcal{D}}
\newcommand{\cE}{\mathcal{E}}
\newcommand{\cF}{\mathcal{F}}
\newcommand{\cG}{\mathcal{G}}

\newcommand{\cI}{\mathcal{I}}

\newcommand{\cK}{\mathcal{K}}

\newcommand{\cO}{\mathcal{O}}
\newcommand{\cP}{\mathcal{P}}

\newcommand{\cS}{\mathcal{S}}

\newcommand{\cU}{\mathcal{U}}

\newcommand{\HH}{\mathbf{H}}

\newcommand{\ZZ}{\mathbf{Z}}

\newcommand{\xx}{\mathbf{x}}
\newcommand{\ee}{\mathbf{e}}
\newcommand{\ff}{\mathbf{f}}
\newcommand{\kk}{\mathbf{k}}

\newcommand{\ggamma}{{\boldsymbol{\gamma}}}

\newcommand{\ess}{\mathrm{ess}}

\newcommand{\CZZ}{\check \mathbf{Z}}
\newcommand{\CZ}{\check Z}

\newcommand{\Cpty}{\mathrm{Cap}}

\theoremstyle{plain} \newtheorem{Theo}{Theorem}[section]
\theoremstyle{plain} \newtheorem{Lemma}[Theo]{Lemma}
\theoremstyle{plain} \newtheorem{Cor}[Theo]{Corollary}
\theoremstyle{plain} \newtheorem{Prop}[Theo]{Proposition}
\theoremstyle{remark} \newtheorem{Rem}[Theo]{Remark}
\theoremstyle{definition} 
\theoremstyle{definition} \newtheorem{Def}[Theo]{Definition}
\theoremstyle{remark} \newtheorem{Not}[Theo]{Notation}
\theoremstyle{definition} \newtheorem{Ass}[Theo]{Assumption}

\begin{document}

\title[Continuum--sites stepping--stone models]{Continuum--sites stepping--stone models, coalescing exchangeable partitions, and random trees}

\author[P. Donnelly]{Peter Donnelly}
\email{donnelly@stats.ox.ac.uk}
\thanks{Donnelly supported in part by NSF grant DMS-9505129,
and UK EPSRC grants B/AF/1255 and GR/M14197}

\address{Department of Statistics\\
Oxford University\\
1 South Parks Road\\
Oxford OX1 3TG\\
U.K.}

\author[S.N. Evans]{Steven N. Evans}
\email{evans@stat.berkeley.edu}
\thanks{Evans supported in part by NSF grant DMS-9703845}

\address{Department of Statistics \#3860 \\
 University of California at Berkeley \\
367 Evans Hall \\
Berkeley, CA 94720-3860 \\
U.S.A.}

\author[K. Fleischmann]{Klaus Fleischmann}
\email{fleischmann@wias-berlin.de}

\address{Weierstrass Institute for Applied Analysis and Stochastics\\
Mohrenstr. 39\\ 
D-10117 Berlin\\
GERMANY}

\author[T.G. Kurtz]{Thomas G. Kurtz}
\email{kurtz@math.wisc.edu}
\thanks{Kurtz supported in part by NSF grant DMS-9626116}

\address{Departments of Mathematics and Statistics\\
University of Wisconsin--Madison\\
Van Vleck Hall\\
480 Lincoln Drive\\
Madison, WI 53706-1388 \\
U.S.A.}

\author[X. Zhou]{Xiaowen Zhou}
\email{xzhou@stat.berkeley.edu}
\thanks{Zhou supported in part by NSF grant DMS-9703845
and a Michel and Lin\'e Lo\`eve Fellowship}

\address{Department of Statistics \#3860 \\
 University of California at Berkeley \\
367 Evans Hall \\
Berkeley, CA 94720-3860 \\
U.S.A.}



\begin{abstract} 
Analogues of stepping--stone models are considered where
the site--space is continuous, the migration process is
a general Markov process, and the type--space is infinite.
Such processes were defined in previous work of the second
author by specifying a Feller transition semigroup
in terms of expectations of suitable functionals
for systems of coalescing Markov processes.
An alternative representation is
obtained here in terms of a limit of
interacting particle systems. 
It is shown that, under  a mild
condition on the migration process, the 
continuum--sites stepping--stone process has
continuous sample 
paths.  The case when the migration
process is Brownian motion on the circle is examined
in detail using a duality relation between coalescing and annihilating
Brownian motion.  This duality relation is also used to show
that a random compact metric space that is naturally associated
to an infinite family of coalescing Brownian motions 
on the circle has
Hausdorff and packing dimension both almost surely equal to $\frac{1}{2}$
and, moreover, this space is capacity equivalent to the 
middle--$\frac{1}{2}$ Cantor set (and hence also to the Brownian
zero set).
\end{abstract}

\maketitle

\section{Introduction}
\label{intro}

{\em Stepping--stone models} originally arose in population genetics.
The simplest version can be described as follows.
There is a a finite or countable collection of sites
(the {\em site--space}).
At each site there is a finite population.
Each population is composed of
individuals who can
be one of two possible genetic types, say A or B.  
At each site
the genetic composition of the population evolves via 
a continuous--time {\em resampling} mechanism.
Independently of each other,
individuals {\em migrate} from one site to another according to
a continuous--time
Markov chain (the {\em migration chain}) on the site--space.

If the number of individuals at each site becomes large, then,
under appropriate conditions, the process describing the proportion of
individuals of type A at the various sites converges to
a diffusion limit.
This limit can be thought of informally as an ensemble of
{\em Fisher--Wright diffusions} 
(one diffusion at each site) that
are coupled together with a  drift determined by the jump
rates of the migration chain (see, for example,
\cite{Shi80}).

A natural refinement of this two--type diffusion model, considered
in \cite{Han90, DaGrVa95}, is 
the corresponding 
{\em infinitely--many--types} model.  Here the Fisher--Wright
processes at each site are replaced by 
mutationless {\em Fleming--Viot processes}
of evolving random
probability measures on a suitable uncountable {\em type--space}
(typically the unit interval $[0,1]$).

Much of the research on such interacting
Fisher--Wright and Fleming--Viot diffusion models 
(see, for example, \cite{BrCoGr86, FlGr94, FlGr96, Kle95})
has centred on their
{\em clustering} behaviour in the case when the space of sites
is either the integer lattice $\bZ^d$ or a 
discrete hierarchical group and
the migration chain is a random walk. That is,
one asks how regions where ``most of the populations are mostly
of one type'' grow and interact with each other.  The primary
tool for analysing this behaviour is the {\em duality} (in the
sense of duality of martingale problems) between these models
and sytems of delayed coalescing random walks that was first
exploited by \cite{Shi80}.

One of the factors that lead to interesting clustering is
the scaling behaviour of the migration process.  However, because
random walks on $\bZ^d$ in the
domain of attraction of a stable law
and their analogues on discrete hierarchical group
only have approximate scaling, the role that scaling plays
is somewhat obscured.  In order to make the
effect of scaling clearer, related two--type
models were considered in \cite{EvFl96} in the hierarchical group setting.
In essence, the processes in \cite{EvFl96}
are the result of taking
a further limit in which one ``stands back'' from the
site--space so that the
discrete hierarchical group approaches a continuous one and the
random walk converges to a ``stable'' L\'evy process
on the continuous hierarchical group that does have exact rescaling.
These {\em continuum--sites, two--type stepping--stone models}
have as their state--space the collection of 
measurable functions $x$ from
the site--space (that is, the continuous hierarchical group)
into $[0,1]$.  For a state $x$ and site $e$, the value
$x(e)$ is 
interpreted intuitively as the proportion of the population
at the site that is of type A.  

One of the noteworthy 
feature of \cite{EvFl96} is that the limit models are defined by
specifying moment--like quantities for the associated Feller
transition semigroup in terms of systems of {\em delayed or
instantaneously coalescing L\'evy processes}, using formulae
that are analogues of the duality relations between the discrete--sites 
models and delayed coalescing random walks mentioned above
(see Theorems 3 and 4 of \cite{EvFl96}).
In particular, the limit models are {\bf not} defined 
infinitesimally via a generator,
SDE/SPDE, or martingale problem formulation analogous to that
of the discrete--sites models.  We note, however, that 
it should be possible to ``stand back''
in a similar manner from a discrete--sites
model where the migration chain
is simple random walk on $\bZ$ 
and obtain the process considered in \cite{MuTr95}:
this process is
constructed there as an SPDE on $\bR$ but is also dual to
delayed coalescing Brownian motions via the same sort of formulae
considered in \cite{EvFl96}.  However, the
processes in \cite{EvFl96} that have their semigroups defined
in terms of instantaneously coalescing L\'evy processes do not
appear to have even a very informal interpretation as SPDE--like
objects.  Rather, a typical value for such a process is a
function $x$ such that $x(e) \in \{0,1\}$ for all sites
$e$, and so such processes are more like continuum analogues
of particle systems (see Theorem 6 of \cite{EvFl96}).

The programme of defining continuum--sites models
in terms of ``duality'' formulae using instantaneously
coalescing Markov processes was continued in \cite{Eva97}
(see Section \ref{defsteppingstone} below for a recapitulation).  
There the infinitely--many--types case was considered
and the migration process (that is, the process used to
build the coalescing system) was taken to be a general
{\em Borel right process} subject only to a duality condition
(with duality here taken in the sense of the general
theory of Markov processes).  Now a state of
the process, which we denote from
now on by $X$, will be a function $x$ from the
site--space into the collection of probability measures on an
uncountable type--space.  For a state $x$, a site $e$ and
a measurable subset $G$ of the type--space, the
value $x(e)(G)$ is interpreted intuitively
as the proportion of the population
at the site possessing types from $G$. 

We give a more concrete description of the
infinitely--many--types, continuum--sites process 
$X$ in Section \ref{partconstructX}.
Under suitable conditions on the migration process, we show
(at the level of convergence of finite--dimensional distributions)
that $X$ is the {\em high--density limit} of a family of 
{\em particle systems}
with the following description.
The particles move about in the Cartesian product of the site--space
and the type--space.  The site--space--valued components
of the particles evolve according to independent copies of the
migration process.  When particles collide in the site--space
a type is chosen at random from the types of the particles
participating in the collision and the types of all the participating
particles are changed to this randomly chosen type.

One of the open problems left by \cite{Eva97} was
to determine conditions on the migration  under which the 
process $X$
(which is again a Feller process)
has {\em continuous} rather than just c\`adl\`ag sample paths.
In Theorem 
\ref{ctspathsthm},
we establish the sufficiency of a mild condition to the effect
that the coalescing system doesn't coalesce too rapidly.
The condition holds, for example, for all L\'evy processes on $\bR$.

By the same argument as in the proof of Proposition 5.1 of \cite{Eva97},
it is possible to show that if the migration process is a stable
process on the {\em circle} $\bT$ that hits points, then, for fixed $t>0$, there
almost surely exists a random countable subset $\{k_1, k_2, \ldots\}$
of the type--space such that for Lebesgue almost all $e \in \bT$ 
the probability measure $X_t(e)$ is a point mass at one of the $k_i$.
That is, rather loosely speaking, at each site all individuals
in the population have the same type and the total number of types seen
across all sites is countable.  We improve this result in
Theorem \ref{XlivesinXio} 
for the case of {\em Brownian motion migration} on $\bT$ by showing
that the total number of types is, in fact, 
almost surely finite and such
a result holds simultaneously at all
positive times rather than just for fixed times.

The primary tool used in the proof of Theorem \ref{XlivesinXio}
 is a {\em duality} relation 
between systems of {\em coalescing} and
{\em annihilating Brownian motions}
that is developed in Section \ref{coalannsect}. 
This relation enables us to perform detailed
computations with the system of coalescing Brownian motions that
begins with countably many particles independently and uniformly
distributed on $\bT$.  

The latter process is an interesting
object in its own right.  In particular, it can be used to
define a {\em random metric} on the positive
integers by declaring that the distance
between $i$ and $j$ is the time until the descendents of the
$i^{\rm{th}}$ and $j^{\rm{th}}$ particles at time zero coalesce.
In Theorem \ref{propsofF}
we adapt the methods of \cite{Eva98} to show that the 
completion of the integers in this metric is
almost surely {\em compact}, with {\em Hausdorff} and 
{\em packing dimensions}
both equal to $\frac{1}{2}$.
Moreover, this space is {\em capacity equivalent} to the 
middle--$\frac{1}{2}$ {\em Cantor set} (and hence also to the Brownian
zero set).

\begin{Not}
Write $\bN:=\{1,2,\ldots\}$.
We will adopt the convention throughout that the
infimum of a subset of $\bR$ or $\bN$ is defined to be
$\infty$ when the subset is empty.
\end{Not}

\section{Coalescing Markov processes and labelled partitions}
\label{coalMarkproc}

Suppose that $E$ is a {\em Lusin space} and that $(Z, P^z)$ is 
a {\em Borel right process} on $E$ with semigroup $\{P_t\}_{t \ge 0}$
satisfying $P_t 1 = 1$, $t \ge 0$, so that $Z$ has infinite lifetime
(see \cite{Sha88} for a discussion of Lusin
spaces and Borel right processes).
Suppose that there is another Borel right
process $(\hat Z, \hat P^z)$ with semigroup $\{\hat P_t\}_{t \ge 0}$ and a
diffuse, Radon measure $m \ne 0$ on $(E, \cE)$ such that
for all  non-negative Borel functions on $f,g$ on $E$
we have $\int m(de) \, P_t f(e) g(e) = \int m(de) \, f(e) \hat P_t g(e)$
(our definition of Radon measure is that
given in Section III.46 of \cite{DeMe78}).
The space $E$ is the {\em site--space} and $\hat Z$ is the
{\em migration process}
for the continuum--sites stepping--stone model $X$, 
whereas $Z$ will serve as the basic motion
in the coalescing systems ``dual'' to $X$.

We remark that our assumption on 
the Markov processes
$Z$ and $\hat Z$ is not quite the
usual notion of {\em weak duality} with respect to $m$
(see, for example, Section 9 of \cite{GeSh84});
in order for weak duality to hold we would also require that  
$P^m$--a.s. (resp. $\hat P^m$--a.s.)
the left--limit $Z(t-)$ (resp. $\hat Z(t-)$) exists for all
$t>0$.

The following notation will be convenient for us.  Given
$\ee = (e_1, \ldots, e_n) \in E^n$, for some $n \in \bN$, with $e_i \ne e_j$
for $i \ne j$, let 
$\ZZ^\ee = (Z^{e_1}, \ldots, Z^{e_n})$ be an $E^n$--valued
process defined on some probability space
$(\Omega, \cF, \bP)$ such that $Z^{e_i}$ has the distribution
of $Z$ under $P^{e_i}$ and $Z^{e_1}, \ldots, Z^{e_n}$ are
independent.

We now define the {\em system  of coalescing
Markov processes} $\CZZ^\ee$ associated with $\ZZ^\ee$.
Adjoin a point $\dag$,
the {\em cemetery}, to $E$ to form $E^\dag := E \cup \{\dag\}$.
Construct the
$(E^\dag)^n$--valued process 
$\CZZ^\ee = (\CZ_1^\ee, \ldots, \CZ_n^\ee)$ 
inductively as follows.
Suppose that times $0 =: \tau_0 \le  \ldots \le \tau_k \le \infty$
and sets 
$\{1, \ldots, n\} 
=: \Theta_0 \supseteq \ldots \supseteq \Theta_k \supseteq \{1\}$
have already been defined and that
$\CZZ^\ee$ has already been defined on $[0,\tau_k[$.
If $\tau_k = \infty$, then just set $\tau_{k+1} := \infty$
and $\Theta_{k+1} := \Theta_k$.
Otherwise, put
\begin{equation}
\tau_{k+1} := \inf\{t>\tau_k : \exists i,j \in \Theta_k,
\, i \ne j, \, Z^{e_i}(t) = Z^{e_j}(t)\},
\end{equation}
\begin{equation}
\Theta_{k+1} := 
         \begin{cases}
         \Theta_k,& \text{if $\tau_{k+1} = \infty$},\\
         \left\{i \in \Theta_k : \not\exists j < i, \, j \in \Theta_k, \, 
           Z^{e_i}(\tau_{k+1}) = Z^{e_j}(\tau_{k+1}) \right\},
         & \text{otherwise},
         \end{cases}
\end{equation}
and 
\begin{equation}
\CZ_i^\ee(t) := 
       \begin{cases}
       Z^{e_i}(t), \, \tau_k \le t < \tau_{k+1}, & \text{if $i \in \Theta_k$},\\
       \dag, \, \tau_k \le t < \tau_{k+1}, & \text{otherwise}.
       \end{cases}
\end{equation}
In other words, the coordinate processes of 
the coalescing Markov process $\CZZ^\ee$ evolve as
independent copies of $Z$ until they collide.  
When two or more coordinate processes collide
(which happens at one of the times $\tau_\ell$ with $0 < \tau_\ell < \infty$),
the one with the smallest index ``lives on'' while the other
coordinates involved in the collision are sent to the cemetery $\dag$.
The set $\Theta_k$ is the set of coordinates that are
still alive at time $\tau_k$.  
As the following lemma shows,
for $m^{\otimes n}$--a.e. $\ee$ almost surely
only one coordinate process of  $\ZZ^\ee$ is sent to the cemetery at
a time in the construction of
$\CZZ^\ee$.  (Recall that  $(Z, P^z)$ is said to be
a {\em Hunt process} if $Z$ has
c\`adl\`ag sample paths and is also {\em quasi--left--continuous};
that is, if 
whenever $T_1 \le T_2 \le \ldots$ are stopping times for $Z$
and $T = \sup_n T_n$, then $P^z\{\lim_n Z(T_n) = Z(T), \, T < \infty\}
= P^z\{T < \infty\}$ for all $z \in E$.)

\begin{Lemma}
\label{triple}
Let $Y$ be an $E$--valued Markov process on some probability
space $(\bar \Omega, \bar \cF, \bar \bP)$ with 
the same law as $Z$ under $P^q := \int_E q(dz) \, P^z$,
where $q$ is a probability measure on $(E, \cE)$ that is
absolutely continuous with respect to $m$.
Let $(T,V)$ be a $[0,\infty[ \times E$--valued random variable
that is independent of $Y$.  Then
$\bar \bP\{Y(T) = V\}  = 0$.
Moreover, if $Z$ is Hunt process, then
$\bar \bP\{Y(T-) = V\} = 0$, also.
A similar result holds with $Y$ replaced by a process
$\hat Y$ with the same law 
as $\hat Z$ under $\hat P^q := \int_E q(dz) \, \hat P^z$.
\end{Lemma}

\begin{proof}
For fixed $t \ge 0$ and $v \in E$ we have,
writing $h$ for the Radon--Nikodym derivative of
$q$ with respect to $m$,
\begin{equation}
\bar \bP\{Y(t) = v\} = \int_E m(dz) \, h(z) P_t \bone_{\{v\}} (z)
= \int_E m(dz) \, \bone_{\{v\}} (z) \hat P_t h(z)
= 0
\end{equation}
by the duality assumption and the assumption that $m$ is diffuse.
Moreover, under the Hunt assumption, 
\begin{equation}
\bar \bP\{Y(t) \ne Y(t-)\} = 0.
\end{equation}
The result now follows by Fubini.
\end{proof}

It will be convenient to embellish $\CZZ^\ee$ somewhat
and consider an enriched process $\zeta^\ee$ defined below
that keeps track of which particles have collided with each
other.

Let $\Pi_n$ denote the set of {\em partitions} of
$\bN_n := \{1, \ldots, n\}$.  That is, an element 
$\pi$ of $\Pi_n$ is a collection 
$\pi = \{A_1, \ldots, A_h\}$ of subsets of $\bN_n$
with the property that 
$\bigcup_i A_i = \bN_n$
and $A_i \cap A_j = \emptyset$ for $i \ne j$.  
The sets $A_1, \ldots A_h$ are the {\em blocks} of the partition $\pi$.
Equivalently,
we can think of $\Pi_n$ as the set of equivalence relations on $\bN_n$
and write $i \sim_\pi j$ if $i$ and $j$ belong to the
same block of $\pi \in \Pi_n$.

An {\em $E$--labelled partition} of
$\bN_n$ is a collection
\begin{equation}
\lambda = \{(A_1,e_{A_1}), \ldots, (A_h,e_{A_h})\}, 
\end{equation}
with
$\{A_1, \ldots, A_h\} \in \Pi_n$, 
$\{e_{A_1}, \ldots, e_{A_h}\} \subseteq  E$,
 and $e_{A_i} \ne e_{A_j}$ for $i \ne j$.
Let $\Lambda_n$ denote the set of $E$--labelled partitions of $\bN_n$.
Put $\alpha(\lambda) := \{A_1, \ldots, A_h\}$ and
$\varepsilon(\lambda) := (e_A)_{A \in \alpha(\lambda)}$.

For $\ee \in E^n$ with $e_i \ne e_j$
for $i \ne j$, we wish to define a $\Lambda_n$--valued 
process $\zeta^\ee$ 
(the {\em process of coalescing Markov labelled partitions}) with
the following intuitive description.  
The initial value of $\zeta^\ee$ is the labelled
partition $\{(\{1\}, e_1), \ldots, (\{n\}, e_n)\}$.
As $t$ increases, the corresponding partition
$\alpha(\zeta^\ee(t))$ remains unchanged and 
the labels $\varepsilon(\zeta^\ee(t))$
evolve as a vector of independent
copies of $Z$  until immediately
before two (or more) such labels coincide.  At the
time of such a collision, the blocks
of the partition corresponding to the coincident labels are
merged into one block
(that is, they {\em coalesce}). This new block
is labelled with the common element of $E$.
The  evolution then continues in the same way.  

More formally, we will take 
$\zeta^\ee$ to be defined in
terms of $\ZZ^\ee$ as follows
(using the ingredients $\tau_k$ and
$\Theta_k$ that went into the
definition of $\CZZ^\ee$).
The corresponding partition--valued
process $\xi^\ee := \alpha(\zeta^\ee)$
is constant on intervals of the form
$[\tau_k, \tau_{k+1}[$
and $\xi^\ee(\tau_0) := \{\{1\}, \ldots, \{n\}\}$.
Suppose for $k \ge 0$ that 
$\xi^\ee(\tau_0), \ldots, \xi^\ee(\tau_k)$
have been defined and $\tau_{k+1} < \infty$.  Let
$\xi^\ee(\tau_{k+1})$ 
be the partition that is obtained by merging
for each $i \in \Theta_{k+1}$ those blocks of
$\xi^\ee(\tau_k)$ whose least elements $j$ are such that
$Z^{e_i}(\tau_{k+1}) = Z^{e_j}(\tau_{k+1})$.
Thus each block $A$
of $\xi^\ee(\tau_{k+1})$ is such that
the least element $\min A$ of $A$ is the unique element
$i \in A$ for which $\CZ_i^\ee(\tau_{k+1}) \ne \dag$.
The definition of $\zeta^\ee$ is completed by labelling
each block $A$ of the partition 
$\xi^\ee(t)$ with $\CZ_{\min A}^\ee (t) = Z^{e_{\min A}}(t)$.

For $1 \le i \le n$, put $\ggamma^\ee = (\gamma_1^\ee, \ldots, \gamma_n^\ee)$,
where 
\begin{equation}
\gamma_i^\ee(t) := \min\{j : j \sim_{\xi^\ee(t)} i\},
\end{equation}
and write
\begin{equation}
\Gamma^\ee(t) := \{\gamma_i^\ee(t): 1 \le i \le n\}
= \{j : \CZ_j^\ee(t) \ne \dag\}
\end{equation}
for the set of surviving indices at time $t$.
Note that $\Gamma^\ee(\tau_k) = \Theta_k$.

\section{The state-space $\Xi$ of the stepping--stone process}

We need some elementary ideas from the theory of vector measures.
A good reference is \cite{DiUh77}.
Recall the measure space $(E, \cE, m)$
introduced in Section \ref{coalMarkproc},
and let $B$ be a Banach space with norm $\| \cdot \|$.
 We say that a function $\phi: E \rightarrow B$ is {\em simple} if
$\phi= \sum_{i=1}^k x_i 1_{E_i}$ for $x_1, \ldots, x_k \in B$ and
$E_1, \ldots, E_k \in \cE$
for some $k \in \bN$.  We say that a function
$\phi: E \rightarrow B$ is $m$--{\em measurable} if there exists a
sequence $\{\phi_n\}_{n \in \bN}$ of simple functions such that
$\lim_{n \rightarrow \infty} \|\phi_n(e)- \phi(e)\| = 0$ for $m$-a.e.
$e \in E$.

Write $K$ for the  compact, metrisable coin-tossing space
$\{0,1\}^{\bN}$ equipped with
the  product topology, and let $\cK$ denote the corresponding
Borel $\sigma$-field.  Equivalently, $\cK$ is the $\sigma$-field
generated by the cylinder sets.

Write $M(K)$ for the Banach space of finite signed measures on $(K,\cK)$
equipped with the total variation norm $\| \cdot \|_{M(K)}$. 
Let $L^\infty(m, M(K))$ denote the space of (equivalence classes
of) $m$-measurable maps $\mu: E \rightarrow M(K)$
such that $\ess \sup \{\|\mu(e)\|_{M(K)} : e \in E\} < \infty$,
and equip  $L^\infty(m, M(K))$ with
the obvious norm to make it a Banach space.

Write $C(K)$ for the Banach space of continuous functions on
$K$ equipped with the usual supremum norm $\|\cdot\|_{C(K)}$.
Let $L^1(m, C(K))$, denote the Banach space of (equivalence classes
of) $m$-measurable  maps $\mu: E \rightarrow C(K)$
such that $\int m(de) \, \|\mu(e)\|_{C(K)} < \infty$, and equip
$L^1(m, C(K))$ with
the obvious norm to make it a Banach space.

{F}rom the discussion at the beginning of \S IV.1 in \cite{DiUh77}
and the fact that $M(K)$ is isometric to the
dual space of $C(K)$
under the pairing
$(\nu, y) \mapsto \langle \nu, y \rangle = \int \nu(dk) \, y(k)$, 
$\nu \in M(K)$, $y \in C(K)$,
we see that $L^\infty(m, M(K))$ is isometric to
a closed subspace of the dual of $L^1(m, C(K))$
under the pairing
$(\mu,x) \mapsto \int m(de) \, \langle \mu(e) , x(e) \rangle$,
$\mu \in L^\infty(m,M(K))$, $x \in L^1(m,C(K))$.
Write $M_1(K)$ for the closed
subset of $M(K)$ consisting of probability measures,
and let $\Xi$ denote the closed subset of $L^\infty(m, M(K))$
consisting of (equivalence classes of) maps with values in $M_1(K)$.
{F}rom Corollary V.4.3 and Theorem V.5.1 of \cite{DuSc58}
we see that, as $L^1(m, C(K))$ is separable,
$\Xi$ equipped with the relative weak$^*$ topology
is a compact, metrisable space.  From now on,
we always take $\Xi$ to be equipped with 
the relative weak$^*$ topology. 

We think of the set $K$ as the space of
possible {\em types} in the infinitely--many--types,
continuum--sites, stepping--stone model $X$
we will define in  Section \ref{defsteppingstone}.
As we remarked in Section \ref{intro},
the type--space for infinitely--many--types models
is usually taken to be $[0,1]$.  However, from a
modelling perspective any uncountable set is equally suitable,
and, as pointed out in \cite{Eva97}, the set $K$
is technically easier to work with.
The set $E$ is the corresponding space of {\em sites}.
The intuitive interpretation is
that $\mu \in \Xi$ describes an ensemble of populations
at the various sites: $\mu(e)(L)$ is the ``proportion of
the population at site $e \in E$ that has
a type belonging to the set $L \in \cK$''.

\begin{Rem}
\label{samemeasstruct}
One can think of $\Xi$ as a subset of the
space of Radon measures on $E \times K$ by identifying
$\mu \in \Xi$ with the measure that assigns mass
$\int_A m(de) \, \mu(e)(B)$ to the set $A \times B$,
where $A \in \cE$ and $B \in \cK$.  The topology we are using on
$\Xi$ is not the same as the trace of the usual topology of
vague convergence of Radon measures.  However, the corresponding
Borel $\sigma$--fields do coincide.  In particular, we can
think of $\Xi$--valued random variables as random Radon measures
on $E \times K$.
\end{Rem}

For $n \in \bN$ let $M(K^n)$ (respectively, $C(K^n)$)
denote the Banach space
of finite signed measures (respectively, continuous functions)
on the Cartesian product
$K^n$ with the usual norm 
$\|\cdot\|_{M(K^n)}$ (respectively, $\|\cdot\|_{C(K^n)}$).
With a slight abuse of notation, write 
$\langle \, \cdot \,, \, \cdot \,\rangle$
for the pairing between these two spaces.

\begin{Def}
\label{defI}
Given  $\phi \in L^1(m^{\otimes n}, C(K^n))$,
define $I_n(\, \cdot \,; \phi) \in C(\Xi)$ 
(:= the space of continuous real--valued functions on $\Xi$) by 
\begin{equation}
\begin{split}
I_n(\mu;\phi) & := \int_{E^n} m^{\otimes n}(d\ee) \, 
\langle \bigotimes_{i=1}^n \mu(e_i),  \phi(\ee) \rangle \\
& = \int_{E^n} m^{\otimes n}(d\ee) \, \int_{K^n} \bigotimes_{i=1}^n
\mu(e_i)(dk_i) \phi(\ee)(\kk), \quad \mu \in \Xi.\\
\end{split}
\end{equation}
Write $I$ for $I_1$.
\end{Def}

\section{Definition of the stepping--stone process $X$}
\label{defsteppingstone}

Theorem \ref{defX} below is Theorem 4.1 of \cite{Eva97}.  
As discussed in Section \ref{intro}, it is motivated by the
characterisation of infinitely--many--types, discrete--sites
stepping--stone processes  via duality with systems of
delayed coalescing continuous--time Markov chains (see \cite{DaGrVa95}
and \cite{Han90}). Recall that $(\Omega, \cF, \bP)$ is
the probability space on which the processes $\ZZ^\ee$, $\CZZ^\ee$,
$\zeta^\ee$, $\xi^\ee$, {\em et cetera} are defined.

\begin{Theo}
\label{defX}
There exists
a unique, Feller, Markov semigroup 
$\{Q_t\}_{t \ge 0}$ on $\Xi$ such that
for all $t \ge 0$, $\mu \in \Xi$,
$\phi \in   L^1(m^{\otimes n}, C(K^n))$, $n \in \bN$,
we have 
\begin{equation}
\label{defQt}
\begin{split}
& \int Q_t(\mu, d\nu) I_n(\nu;\phi) \\
& \quad = \int_{E^n} m^{\otimes n}(d\ee)
\bP
\left[
\int \bigotimes_{j \in \Gamma^\ee(t)} \mu(Z_j^\ee(t))(dk_j)
\phi(\ee)(k_{\gamma_1^\ee(t)}, \ldots, k_{\gamma_n^\ee(t)})
\right]. \\
\end{split}
\end{equation}
Consequently, there is a Hunt process, $(X, \bQ^\mu)$, with 
state-space $\Xi$ and transition semigroup
$\{Q_t\}_{t \ge 0}$.
\end{Theo}

\begin{Rem}
The integrand $\bP[\cdots]$ in (\ref{defQt}) should be interpreted
as $0$ on the $m^{\otimes n}$--null set
of $\ee$ such that $e_i = e_j$ for some pair
$(i,j)$.  The integral inside the $[ \; ]$ is over
a Cartesian
product of copies of $K$, with the copies indexed by the
elements of $\Gamma^\ee(t)$.
\end{Rem}

\begin{Rem} 
\label{ZndefX}
The following equivalent formulation
of Theorem \ref{defX} will be useful.  For $n \in \bN$ let 
$\ZZ^{[n]} = (Z_1^{[n]}, \ldots, Z_n^{[n]})$ 
be an $E^n$--valued process defined
on a $\sigma$-finite measure space $(\Omega^{[n]}, \cF^{[n]}, \bP^{[n]})$,
with 
\begin{equation}
\bP^{[n]}\{\ZZ^{[n]} \in A\} := \int m^{\otimes n}(d\ee) \bP\{\ZZ^\ee \in A\}.
\end{equation}
Define $\CZZ^{[n]}$, $\xi^{[n]}$,
$\ggamma^{[n]}$ and $\Gamma^{[n]}$ from $\ZZ^{[n]}$ in the same
manner that
$\CZZ^\ee$, $\xi^\ee$, $\ggamma^\ee$ and $\Gamma^\ee$ 
were defined from $\ZZ^\ee$.  The right--hand side of
(\ref{defQt}) is just
\begin{equation}
\bP^{[n]}
\left[
\int \bigotimes_{j \in \Gamma^{[n]}(t)} \mu(Z_j^{[n]}(t))(dk_j)
\phi(\ZZ^{[n]}(0))(k_{\gamma_1^{[n]}(t)}, \ldots, k_{\gamma_n^{[n]}(t)})
\right]. 
\end{equation}
\end{Rem}

\begin{Rem}
\label{extnofXdef}
As we noted in Remark \ref{samemeasstruct}, we can think of
the process $X$ as taking values in the space of Radon measure
on $E \times K$ by identifying $X_t$ with the random measure
that assigns mass $\int_A m(de) X_t(e)(B)$ to the set
$A \times B$, where $A \in \cE$ and $B \in \cK$.
A standard monotone class argument shows that if
$\psi$ is any non-negative Borel function on
$E^n \times K^n$, then 
\begin{equation}
\begin{split}
&\bQ^\mu\left[
\int_{E^n} m^{\otimes n}(d\ee) \, 
\int_{K^n} \bigotimes_{i=1}^n X_t(e_i)(dk_i) \, \psi(\ee,\kk)\right] \\
& \quad =
\int_{E^n} m^{\otimes n}(d\ee)
\bP
\left[
\int \bigotimes_{j \in \Gamma^\ee(t)} \mu(Z_j^\ee(t))(dk_j)
\psi(\ee, k_{\gamma_1^\ee(t)}, \ldots, k_{\gamma_n^\ee(t)})
\right] \\
& \quad =
\bP^{[n]}
\left[
\int \bigotimes_{j \in \Gamma^{[n]}(t)} \mu(Z_j^{[n]}(t))(dk_j)
\psi(\ZZ^{[n]}(0), k_{\gamma_1^{[n]}(t)}, \ldots, k_{\gamma_n^{[n]}(t)})
\right]. \\
\end{split}
\end{equation}
\end{Rem}

\section{A particle construction for the stepping--stone model $X$}
\label{partconstructX}

In this section we first construct a finite particle 
model in which particles
move through $E \times K$, 
where we recall that $E$ is our {\em site--space}
and $K$ is our {\em type--space}.
The $E$--valued components of the particles move independently
according to the dynamics of the migration process $\hat Z$.  
The particles interact only when they are 
located at the same site in $E$, and the interaction that 
occurs is that the type of one of the particles is 
replaced by the type of the other.  The particle whose 
type ``wins'' is chosen  at random from the two 
particles, with both outcomes equally likely. 
For our purposes here,
we assume that the types are constant except for these 
replacement interactions, although we could allow 
``mutation'' of the types between the replacement 
interactions.

Under suitable conditions on the migration process, we then 
pass to a high--density limit and obtain a
process taking values in the space of Radon measures
$\rho$ on $E \times K$ with the property that
$\rho(A \times K) = m(A)$ for $A \in \cE$.
Recalling Remark \ref{samemeasstruct}, we can think of the limit
model as a $\Xi$--valued process, and we establish
that as such it has the same finite--dimensional distributions
as the continuum--sites stepping--stone process $X$.

Throughout this section we will work on a probability space
$(\hat \Omega, \hat \cF, \hat \bP)$ and
we will assume the following hypothesis
(the definition of a Hunt process is recalled in
Section \ref{coalMarkproc}).

\begin{Ass} 
\label{Huntass}
The processes $Z$ and $\hat Z$ are Hunt.
\end{Ass}

For completeness, we recall the following definition 
and some of its consequences.

\begin{Def}
Let $(S,\cS)$ be a measurable space, and let $\nu$ be a 
$\sigma$-finite measure on $\cS$.  
Say that a map $N$ from
$\hat \Omega$ into the
collection of measures on
$(S, \cS)$
is a {\em Poisson random measure }
with {\em mean measure} $\nu$ if

\begin{itemize} 

\item[a)]For each $A \in \cS$,
$N(A)$ is a $\{0,1, \ldots, \infty\}$--valued
random variable.

\item[b)]For each $A\in \cS$ with $\nu (A)<\infty$, 
the random variable $N(A)$ is 
Poisson distributed with parameter $\nu (A)$.

\item[c)]For 
$A_1,A_2,\ldots\in \cS$ disjoint, the random variables
$N(A_1), N(A_2),\ldots$ are independent.

\end{itemize}
\end{Def}

\begin{Rem}
Assume that $\nu$ is diffuse.  Then for $x\in S$, $N(\{x\})$ must be 
zero or one, and so we can identify $N$ with its support.  We 
will write $x \in N$ if $N(\{x\})=1$.  Note that
\begin{equation}
\hat \bP\left[
\int_E
N(dx) f(x) \right]
=
\hat \bP\left[
\sum_{x\in N} f(x) \right]
=
\int_E
\nu (dx) f(x),\;
f \in L^1(\nu ),
\end{equation}
and more generally, 
for $f\in L^1(\nu^{\otimes n})$,
\begin{equation}
\label{momid}
\hat \bP\left[
\sum_{\begin{array}{c}
x_1,\ldots ,x_n\in N\\
x_i \ne x_j, \, i \ne j\end{array}
}
f(x_1,\ldots ,x_n) \right]
=
\int_{E^n} \nu^{\otimes n}(d\xx)
f(x_1,\ldots ,x_n).
\end{equation}
\end{Rem}

\subsection{Finite particle systems}

Fix a non-zero diffuse finite measure $\nu_0$ on $E$
and a probability kernel $\mu: E \times \cK \rightarrow [0,1]$.
Write $D_E[0,\infty[$ for the Skorohod space of
c\`adl\`ag $E$--valued paths and let 
$\widehat M$ denote a Poisson random measure on
$D_E[0,\infty[ \times K$ with mean measure 
\begin{equation}
F\times G \mapsto \int \nu_0(dz) \, \hat P^z(F) \, \mu (z,G)
\end{equation}
(recall that $\hat P^z$ is the law of $\hat Z$ starting at $z \in E$).
Thus the push--forward of $\widehat M$ by the map
$(\zeta, k) \mapsto \zeta(0)$ ($:=$ the value of
the path $\zeta$ at time $0$) is a Poisson random measure on $E$
with mean measure $\nu_0$.  More generally,
the push--forward of $\widehat M$ by the map
$(\zeta, k) \mapsto \zeta(t)$ is a Poisson random measure on $E$
with mean measure $\nu_t$, where 
$\nu_t(H)=\int_E \nu_0(dz) \, \hat {P}_t(z,H)$.
We assume that $\nu_t$ is diffuse for each $t \ge 0$.  
By our duality assumption,
this will certainly be the case if $\nu_0$ is absolutely
continuous with respect to $m$.

Enumerate the atoms of $\widehat M$ as
$(\hat Z_1, \kappa_1^0), \ldots, (\hat Z_J, \kappa_J^0)$
in such a way that the conditional distribution of
this collection given $J=j$ is that of $j$ i.i.d.
$D_E[0,\infty[ \times K$--valued
 random variables with common distribution 
\begin{equation}
F\times G \mapsto \nu_0(E)^{-1} \int \nu_0(dz) \, \hat P^z(F) \, \mu (z,G).
\end{equation}

We wish to define a collection
$\kappa_1, \ldots, \kappa_J$
of $K$--valued processes in such way that the collection
$(\hat Z_1, \kappa_1), \ldots, (\hat Z_J, \kappa_J)$
has the dynamics described above: that is,
we think of $\kappa_i(t)$ as the type of the particle $\hat Z_i$
at time $t$, 
and,  after two or more such particles collide, 
the particles participating in the collision must be of the same 
type with the common type selected at 
random from among the types of the participating particles 
(with each possible outcome equally likely).

Suppose that on the probability space
$(\hat \Omega, \hat \cF, \hat \bP)$
we also have defined for each $k \in \bN$
a collection $\{\theta_{ik}, \, i \in \bN \}$ of i.i.d. 
random variables uniformly distributed on $[0,1]$.
We will implement a specific construction of the $\theta_{ik}$ below.
Define $\kappa_1, \ldots, \kappa_J$
and times $\hat \tau_0 \le \hat \tau_1 \le \ldots$
(with $\hat \tau_k < \hat \tau_{k+1}$ when $\hat \tau_k < \infty$) 
as follows.

Put $\kappa_i(0) = \kappa_i^0$ and $\hat \tau_0 = 0$.
Suppose that $\hat \tau_0, \ldots, \hat \tau_k$ have already
been defined and, for $1 \le i \le J$, 
the processes $\kappa_i$ has been defined
on $[0, \hat \tau_k]$
(or $[0, \infty[$ if $\hat \tau_k = \infty$).
If $\hat \tau_k = \infty$, then the definition of $\kappa_i$,
$1 \le i \le J$, is complete and just define $\hat \tau_\ell = \infty$
for $\ell>k$.
Suppose, then, that $\hat \tau_k < \infty $. Put
\begin{equation}
\hat \tau_{k+1}:=\inf\{t>\hat \tau_k:\hat Z_i(t)=\hat Z_j(t), \,
\kappa_i(\hat \tau_k)\neq\kappa_j(\hat \tau_k), \text{ some }i\neq j\}.
\end{equation}
Put $\kappa_i(t) := \kappa_i(\hat \tau_k)$ for 
$\hat \tau_k \le t < \hat \tau_{k+1}$
and $1 \le i \le J$.
If $\hat \tau_{k+1} = \infty$, then this completes the definition of 
$\kappa_i$, $1 \le i \le J$.  Otherwise, if $\hat \tau_{k+1} < \infty$, 
then define $\kappa_i(\hat \tau_{k+1})$, $1 \le i \le J$,
as follows.
Let $\hat \Gamma_i(\hat \tau_{k+1}) := 
\{j:\hat Z_j(\hat \tau_{k+1})=
\hat Z_i(\hat \tau_{k+1})\}$, and let 
$\hat \gamma_i(\hat \tau_{k+1})\in\hat\Gamma_i(\hat \tau_{k+1})$ 
satisfy  
$\theta_{\hat \gamma_i(\hat \tau_{k+1}),k+1}\leq\theta_{j,k+1}$ 
for all $j\in\hat\Gamma_i(\hat \tau_{k+1})$.  We 
set $\kappa_i(\hat \tau_{k+1})
=\kappa_{\hat \gamma_i(\hat \tau_{k+1})}(\hat \tau_k)$. 
Note that if $\hat Z_i(\hat \tau_{k+1})=\hat Z_j(\hat \tau_{k+1})$, then 
$\hat \Gamma_i(\hat \tau_{k+1})=\hat \Gamma_j(\hat \tau_{k+1})$, 
$\hat \gamma_i(\hat \tau_{k+1})
=\hat \gamma_j(\hat \tau_{k+1})$, 
and $\kappa_i(\hat \tau_{k+1})=\kappa_j(\hat \tau_{k+1})$.

Our requirement that the types of colliding particles be changed to
a type independently and uniformly selected from those of the participants
in a collision will be met if for
$k \in \bN$, the collection $\{\theta_{ik}, \, i \in \bN \}$ is
independent of
$\cF^{\hat Z}_{\hat \tau_k}\vee \cF^{\kappa}_{\hat \tau_{k-1}}$,
where $\{\cF^{\hat Z}_t\}_{t \ge 0}$ is the filtration
generated by  $(\hat Z_1, \ldots \hat Z_J)$ and
$\{\cF^{\kappa}_t\}_{t \ge 0}$ is the filtration generated 
by $(\kappa_1, \ldots, \kappa_J)$.  In particular, 
the distribution of the process $(\hat Z,\kappa )$ will be the same 
regardless of how we define the $\{\theta_{ik}\}$ as long as for 
each $k$, the conditional distribution of $\{\theta_{ik}\}$ given 
$\cF^{\hat Z}_{\hat \tau_k}\vee \cF^{\kappa}_{\hat \tau_{k-1}}$ 
is i.i.d. uniform on $[0,1]$.

We note 
that $\hat \bP$-a.s. 
there exists $\ell \in \bN$ such that $\hat \tau_\ell=\infty$,
so that the above construction does indeed lead to
a value of $\kappa_i(t)$, $1 \le i \le J$,
for all $t \ge 0$.  To see this, let 
$R_h(t)=\{1 \le j \le J:\kappa_j(t)=\kappa_h(t)\}$
for $1 \le h \le J$
and $0 \le t < \sup_k \hat \tau_k$.  Since there are only finitely 
many particles, $\hat \bP\{R_h(\hat \tau_k)\subseteq R_h(\hat \tau_{k+1})\subset\cdots 
\, | \, \cF_{\hat \tau_k}^{\hat Z} 
\vee \cF_{\hat \tau_k}^{\kappa}\}\geq 2^{-J}>0$.  
Consequently, either there exists $\hat \tau_k<\infty$ such that 
$R_h(\hat \tau_k)=\{1,\ldots ,J\}$ or there exists a time after which
$\hat Z_h$ 
does not collide with any particle having a different 
type.  

Now we will give an explicit 
construction of the $\{\theta_{ik}\}$ which leads to a useful 
construction of our particle system $(\hat Z_1, \kappa_1),
\ldots, (\hat Z_J, \kappa_J)$.
We assign to each particle a distinct $[0,1]$--valued 
initial {\em level} $U_i^0$, 
$1 \le i \le J$,
at time $0$ and use these initial levels to define
a family of $[0,1]$--valued processes  of levels $\{U_i(t)\}_{t \ge 0}$,
$1 \le i \le J$.
The $\{\theta_{ik}\}$ will be defined using these level processes.  We 
will assume that the conditional distribution of
$\{U_i^0\}$ given $\widehat M$
is that of $J$
i.i.d. random variables uniformly 
distributed on $[0,1]$.  This assumption implies
\begin{equation}
\sum_{i=1}^J\delta_{(\hat Z_i,\kappa_i^0,U_i^0)}
\end{equation}
is a Poisson random measure with mean measure
\begin{equation} 
\label{meanmeasspec}
F\times G\times H \mapsto \int_E \nu_0(dz) \, \hat P^z(F) \, \mu(z,G) \, l(H),
\end{equation} 
where $l$ denotes 
Lebesgue measure on $[0,1]$.  We define $\theta_{i1}:=U_i^0$.
For $0\leq t<\hat \tau_1$, set $U_i(t):=U_i(0):=U_i^0$. 
If $\hat \tau_1 < \infty$ and $|\hat \Gamma_i(\hat \tau_1)|=1$,
then put $U_i(\hat \tau_1):=U_i(\hat \tau_1-)$.   
If $\hat \tau_1 < \infty$ and $\hat \Gamma_i(\hat \tau_1)=\{i_1,\ldots ,i_n\}$, 
$n>1$, then 
put $U_{i_l}(\hat \tau_1):=U_{i_{\sigma_l}}(\hat \tau_1-)$ where  $\sigma_1,\ldots 
,\sigma_n$ is a uniform random 
permutation of $1,\ldots ,n$ selected independently of all other 
quantities.  Observe that $U_1(\hat \tau_1),\ldots ,U_J(\hat \tau_1)$  are 
conditionally i.i.d. uniform on $[0,1]$ given $\cF^{\hat Z}_{\hat \tau_
2}\vee \cF^{\kappa}_{\hat \tau_1}$. Define 
$\theta_{i2}:=U_i(\hat \tau_1)$.  Put $U_i(t):=U_i(\hat \tau_1)$, $\hat \tau_1<t<\hat \tau_
2$.  We continue 
inductively, at each time $\hat \tau_k < \infty$
randomly permuting the 
levels with indices in  each $\hat \Gamma_i(\hat \tau_k)$ and defining 
$\theta_{i,k+1}=U_i(\hat \tau_k)$.

Although the level assigned to a particle may change at 
the time of a collision, since these changes only involve 
the permuation of the assignment of the levels, the set of 
levels is the fixed random set $\cU:=\{U_i^0\}$.  Consequently, we could 
index the particles and their types by their 
corresponding level; that is, for $u\in \cU$, define
$\hat Z_u(t)=\hat Z_i(t)$ and $\kappa_u(t)=\kappa_i(t)$ if and only if 
$U_i(t)=u$.  
Since the particle assigned to level $u$ changes only when 
the newly assigned particle is at the same location as 
the previously assigned particle, the strong Markov 
property implies that the processes
$\{\hat Z_u, \, u \in \cU\}$ are conditionally
independent given  $\cU$ and $\{\hat Z_u(0), \, u\in \cU\}$,
and conditionally each $\hat Z_u$ is a Markov process with 
transition semigroup $\{\hat {P}_t\}$.  Note that 
\begin{equation} 
\hat \tau_{k+1}=\inf\{t>\hat \tau_k:\hat Z_u(t)=\hat Z_v(t),
\, \kappa_u(t-)\neq\kappa_v(t-),\text{ some }u\neq v\},
\end{equation}
and if we define
$\hat \Gamma_u(\hat \tau_k):=\{v\in \cU:\hat Z_v(\hat \tau_k)=\hat Z_u(\hat \tau_k)\}$ 
and $\hat \gamma_u(\hat \tau_k):=\min(\hat \Gamma_u(\hat \tau_k))$, 
then $\kappa_u(\hat \tau_k)=\kappa_{\hat \gamma_u(\hat \tau_k)}(\hat \tau_k-)$.
That is, if two 
or more particles collide,
the particles involved in the collision 
``look down'' to the lowest level particle at the same 
location, and change types to the type of that particle.
(We note in passing that 
this construction is reminiscent of the ``look down''
construction of the Moran model in \cite{DoKu96}.)
Consequently, if we start with a Poisson random 
measure on $D_E[0,\infty[ \times K \times [0,1]$
\begin{equation}
\sum_{u\in \cU}\delta_{(\hat Z_u,\kappa_u^0,u)}
\end{equation}
with mean measure 
specified by (\ref{meanmeasspec}), then
 a particle model
\begin{equation}
\Psi_t=\sum_{u\in \cU}\delta_{(\hat Z_u(t),\kappa_u(t)
,u)}
\end{equation}
is completely determined by the requirement that 
whenever two or more particles ``collide'' the types of 
the higher level particles involved in the collision 
switches to the type of the lowest level particle in the 
collision.  This observation allows us to extend the 
construction to systems with infinitely many particles 
with mild additional assumptions.

\subsection{Particle systems with stationary location processes}

We now want to extend the construction of the previous
section to arrive at a model in which
the distribution of locations of particles is
stationary in time, and we want to allow for the 
possibility of there being infinitely many particles.
We emphasize that Assumption \ref{Huntass} is in force
throughout this section.

Consider a Poisson random measure
\begin{equation}
\sum_{u\in \cU} \delta_{(\hat Z_u,\kappa_u^0,u)}
\end{equation}
on $D_E[0,\infty[ \times K \times [0,1]$ with 
mean measure 
specified by (\ref{meanmeasspec}) with $\nu_0 = m$.
By the assumption that $m$ is Radon, there exist 
open sets $E_1 \subseteq E_2 \subseteq \cdots$
such that $m(E_n)< \infty$ for all $n \in \bN$
and $E = \bigcup_n E_n$ (of course, if $m(E)$ is finite we
can take $E_n = E$ for all $n \in \bN$).
Put 
$\cU_n=\{u\in \cU:\hat Z_u(0)\in E_n\}$, and note that
\begin{equation}
\sum_{u\in \cU_n}\delta_{(\hat Z_u,\kappa_u^0,u)}
\end{equation}
is a Poisson random measure on 
$D_E[0,\infty[ \times K\times [0,1]$ with finite mean measure 
specified by (\ref{meanmeasspec}) with $\nu_0 = m(\cdot \cap E_n)$.  
As in the previous subsection, we 
can construct a corresponding finite particle model
\begin{equation}
\Psi^{[n]}(t)=\sum_{u\in \cU_n}\delta_{(\hat Z_u(t),\kappa^{[n]}_
u(t),u)}
\end{equation}
and times $\hat \tau_0^{[n]} \le \hat \tau_1^{[n]} \le \ldots$.
We would like to define 
$\Psi_t=\lim_{n\rightarrow\infty}\Psi^{[n]}(t)$;  
however, the type processes $\kappa_u^{[n]}$ 
may not converge without some additional assumptions 
regarding the behavior of the 
migration processes $\{\hat Z_u, \, u \in \cU\}$.

Henceforth, we will also assume the following, which
will ensure that for all $n \in N$ and $t \ge 0$ the expectation
$\hat \bP[ \, |\{u \in \cU : \hat Z_u(s) \in E_n 
\; \text{for some} \; 0 \le s \le t\}| \, ]$ is finite.

\begin{Ass} 
\label{finitecap}
The sequence $\{E_n\}$ of open sets can be chosen
so that 
\begin{equation*}
\hat P^m\{\sigma_{E_n} \le t\} < \infty, 
\; \text{for all $n \in \bN$ and $t > 0$},
\end{equation*}
 where
\begin{equation*}
\sigma_A := \inf\{t \ge 0 : \hat Z(t) \in A\}, \; A \in \cE.
\end{equation*}
\end{Ass}

\begin{Rem}
By our duality assumption, the measure $m$ is stationary
for $\hat Z$.  It follows easily that if for $A \in \cE$
the condition
$\hat P^m\{\sigma_A \le t\} < \infty$ holds for
some $t>0$, then it holds for all $t >0$. Furthermore,
the condition 
$\hat P^m\{\sigma_A \le t\} < \infty$ for
all $t>0$
is also equivalent to
$\hat P^m[\exp(-\lambda \sigma_A)] < \infty$
for all (equivalently, some) $\lambda>0$.
Using this equivalence, the question of whether or
not Asssumption \ref{finitecap} is satisfied becomes
a standard question in capacity theory.
Under our duality assumption and the Assumption \ref{Huntass} that 
$Z, \hat Z$ are Hunt, Assumption \ref{finitecap}
will certainly hold 
(with $\{E_n\}$ any increasing
sequence of relatively compact open sets such that
$\bigcup_n E_n = E$)
if the Lusin space $E$ is locally compact
and $\lambda$-excessive functions for both semigroups
$\{P_t\}$ and $\{\hat P_t\}$ are lower semi--continuous
(see, for example, Remark 2.10 of \cite{Get84}).
In particular, 
Assumption \ref{finitecap} holds if 
$E$ is locally compact and
$Z$ and $\hat Z$
have strong Feller $\lambda$--resolvent operators
(see Exercise II.2.16 of \cite{BlGe68}).   
Also, Assumption \ref{finitecap} holds
when $Z$ and $\hat Z$ are
L\'evy processes on $\bR^d$ and $m$ is Lebesgue measure (see Lemma II.6
of \cite{Ber96}).
\end{Rem}

Fix $t>0$ and $u\in \cU_n$.
Let 
\begin{equation}
\alpha_n(u,t) :=
0\vee\sup\{0<s<t:\hat Z_u(s)=\hat Z_v(s) \text{ some }
v<u, \, v\in \cU_n\},
\end{equation}
and let $\beta_n(u,t)$ be the corresponding value of $v\in \cU_n$, 
with $\beta_n(u,t):=u$ if $\alpha_n(u,t)=0$.   
Lemma \ref{triple}\ 
implies that $\beta_n(u,t)$ is well-defined.  
In general, $\alpha_n(u,t)$ will not be one of the times
$\{\hat \tau_k^{[n]}\}$, but 
we will have
\begin{equation}
\kappa_u^{[n]}(t)=\kappa_u^{[n]}(\alpha_n(u,t))=\kappa_{\beta_n(u,t)}^{[n]}
(\alpha_n(u,t)).
\end{equation}
Define $\beta_{n,u}^t(s):=u$ for $\alpha_n(u,t) <s \leq t$ and 
$\beta_{n,u}^t(s):=\beta_n(u,t)$ for 
$\alpha_n(\beta_n(u,t),\alpha_n(u,t))<s\leq\alpha_n(u,t)$.  This definition extends 
iteratively to determine $\beta_{n,u}^t(s)$ on the interval 
$0\leq s\leq t$ 
with the property that 
\begin{equation}
\kappa_u^{[n]}(t)=\kappa^{[n]}_{\beta_{n,u}^t(s)}(s),
\end{equation}
so, in particular, $\kappa_u^{[n]}(t)=\kappa_{\beta^t_{n,u}(0)}^0$.
  Consequently, 
convergence of $\kappa_u^{[n]}(t)$ is equivalent to convergence of 
$\beta^t_{n,u}$.

For $t>0$ and $u \in \cU$ set 
\begin{equation}
\label{defalphaut}
\alpha (u,t)
:=0\vee\sup\{0<s<t:\hat Z_u(s)=\hat Z_v(s) \text{ some }
v<u, \, v\in \cU\}.
\end{equation}
Let $U$ be a $[0,1]$--valued random variable
that is $\sigma(\cU)$--measurable 
 and takes values in the random set
$\cU$  (that is, $U$ is a
$\sigma(\cU)$--measurable selection from $\cU$). 
By the duality assumption
and the Hunt hypothesis Assumption \ref{Huntass}
(cf. Proposition 15.7 of \cite{GeSh84}),
$t-\alpha (U,t)$ has the same distribution as
\begin{equation} 
\inf\{s>0:Z_U(s)=Z_v(s) \text{ some }v<U, \, v\in \cU\} \wedge t,
\end{equation}
where $\sum_{u\in \cU}\delta_{(Z_u,u)}$ is any Poisson random measure 
with mean measure 
\begin{equation}
F\times H \mapsto \int_E m(dz) \, P^z(F) \, l(H)
\end{equation}
constructed from $\cU$
using suitable further randomisation.
Let $\beta (u,t)$ be the corresponding value of $v\in \cU$
in (\ref{defalphaut}), with 
$\beta (u,t):=u$ if $\alpha (u,t)=0$.  
Define $\beta_u^t(s):=u$ for $\alpha (u,t)<s\;\leq t$ and 
$\beta_u^t(s):=\beta (u,t)$ for 
$\alpha (\beta (u,t),\alpha (u,t))<s\leq\alpha (u,t)$.  Extending this definition 
iteratively, either we determine $\beta_u^t(s)$ on the interval 
$0\leq s\leq t$ and there are only finitely many levels in the 
range of $\beta^t_u$ or there exists $T_u^t \geq 0$ such that 
$\lim_{s\downarrow T_u^t}\beta_u^t(s)=0$.  We show that
this  latter possibility cannot occur.

Suppose that the latter possibility does occur.
As above, let $U$ be a $\sigma(\cU)$--measurable random
variable taking values in $\cU$.
Define $T_U^t$ by analogy with $T_u^t$, with the convention that
$T_U^t := t$ if there are only finitely many levels in
the range of $\beta_U^t$.
Set
\begin{equation}
Z_U^t(r):= \lim_{r' \downarrow r} \hat Z_{\beta_U^t(t-r')}(t-r')
, \; 0 \le r < t - T_U^t.
\end{equation}
Then by the strong Markov property, 
the duality assumption and Assumption \ref{Huntass}, 
$\{Z_U^t(r), \,  0 \le r < t - T_U^t\}$ is a c\`adl\`ag Markov 
process with transition semigroup $\{P_r\}$.  
In particular, the range of this process is almost surely
relatively compact and is contained in one of the $E_n$
for $n$ sufficiently large.
Now, by Assumption \ref{finitecap}, the
cardinality of the set  
\begin{equation}
\{v \in \cU : v < U, \,  \hat Z_v(s) \in E_n \text{ some }
0 \le s \le t\}
\end{equation} 
is $\hat \bP$-a.s. finite for
all $n \in \bN$.  Consequently,
$\hat \bP$-a.s.
there are indeed only finitely many levels in the 
range of $\beta_U^t$
and hence only finitely many levels in
the range of $\beta_u^t$ for all $u \in \cU$.  It follows that 
$\hat \bP$-a.s. we have 
$\lim_{n\rightarrow\infty} \beta^t_{n,u}=\beta^t_u$
 and hence
\begin{equation}
\lim_{n\rightarrow\infty}\kappa_u^{[n]}(t)=\kappa_{\beta_u^t(0)}^0
=: \kappa_u(t)
\end{equation}
for all $u \in \cU$.

If $\beta^t_{u_1}(s)=\beta^t_{u_2}(s)$ for some $0\leq s\leq t$, then 
$\beta^t_{u_1}(s')=\beta^t_{u_2}(s')$ for all $0\leq s'\leq s$.
Moreover, if we define
\begin{equation}
Z^t_u(r)=\lim_{r' \downarrow r}  \hat Z_{\beta^t_u(t-r')}(t-r'),
\; 0 \le r \le t,
\end{equation}
for each $u \in \cU$, then conditional 
on $\cU$ each $Z^t_u$ is a Markov process with transition 
semigroup $\{P_r\}$.  In particular, $\{Z_u^t, \, u\in \cU\}$ form a 
coalescing system of Markov processes, and for $0 \le r \le t$ the 
equivalence relation defined by $u\sim v$ if and only if 
$\beta_u^t(t-r)=\beta_v^t(t-r)$ determines a partition $\{\cU_
k^t(r)\}$ of the 
set of levels $\cU$.  For definiteness, assume that $\cU$ is 
ordered $\cU=\{U_1,U_2,\ldots \}$
where the $U_i$ are $\cU$--measurable random variables, 
and let $\cU_1^t(r)$ be the equivalence 
class containing $U_1$, let $\cU_2^t(r)$ be the equivalence class 
containing the $U_i$ with smallest index not contained in 
$\cU_1^t(r)$, etc.  For each $k$, $Z^t_u(r)$ has the same value for all 
$u\in \cU_k^t(r)$, which we denote by $Z_k^t(r)$.  Then 
$((Z_1^t,\cU_1^t),(Z_2^t,\cU_2^t),\ldots )$ forms a 
coalescing Markov labelled 
partition of $\cU$.  

Since the initial particle types $\{\kappa_u^0\}$ are 
conditionally independent given $\{\hat Z_u, \, u \in \cU\}$ and $\cU$, and 
\begin{equation}
\begin{split}
& \hat \bP\left[g(\kappa_{u_1}^0,\ldots ,\kappa_{u_n}^0) \, 
| \, \{\hat Z_u, \, u \in \cU\}, \, \cU \right] \\
& \quad =\int_{K^n} \mu (\hat Z_{u_1}(0),dk_1)\cdots
\mu (\hat Z_{u_n}(0),dk_n) g(k_1,\ldots ,k_n), \\
\end{split}
\end{equation}
for $u_1,\ldots ,u_n\in \cU$, we have
\begin{equation}
\label{cndid}
\begin{split}
& \hat \bP\left[f(\hat Z_{u_1}(t),\kappa_{u_1}(t),\ldots ,\hat {
Z}_{u_n}(t),\kappa_{u_n}(t)) \, \Big | \, \{\hat Z_u(t)\}, \, \cU \right] \\
& \quad = \HH_t f(\hat Z_{u_
1}(t),\ldots ,\hat Z_{u_n}(t)), \\ 
\end{split}
\end{equation}
where, in the notation of Section \ref{coalMarkproc}, 
\begin{equation}
\label{coalop}
\HH_t f(e_1, \ldots, e_n) :=
\hat \bP
\left[
\int \bigotimes_{j \in \Gamma^\ee(t)} \mu(Z_j^\ee(t), dk_j)
f(e_1, k_{\gamma_1^\ee(t)}, \ldots, e_n, k_{\gamma_n^\ee(t)})
\right]
\end{equation}
for $e_1, \ldots, e_n \in E$ with $e_i \ne e_j$, $i \ne j$.

By (\ref{momid}) and (\ref{cndid}) we have
\begin{equation} 
\label{dualid}
\begin{split}
&\hat \bP\left[
\sum_{
\begin{array}{c}
u_1,\ldots ,u_n\in \cU \\
u_i\neq u_j, \, i \ne j
\end{array}
}
f(\hat Z_{u_1}(t),\kappa_{u_1}(t),\ldots ,\hat Z_{u_n}(t),\kappa_{
u_n}(t)\right]\\
& \quad =\int_{E^n} m^{\otimes n}(d\ee) \, \HH_t f(\ee), \\
\end{split}
\end{equation}
for $f$ a bounded measurable function
on $(E\times K)^n$.  This identity gives a duality
 in the 
sense of (4.4.36) of \cite{EtKu86} between 
the discrete--particle, continuum--sites model and the 
corresponding coalescing Markov labelled partition process.

Write 
\begin{equation}
\Psi^1_t=\sum_{u\in \cU}\delta_{(\hat Z_u(t),\kappa_u(t
),u)}
\end{equation}
and 
\begin{equation}
X^1_t=\sum_{u\in \cU}\delta_{(\hat Z_u(t),\kappa_u(t))}.
\end{equation}
Set $\cF^{X^1}_t=\sigma (X^1_s:s\leq t)$.
The levels of the particles are independent of $\cF_t^{X^1}$, so 
if $f(e,k,u)$ satisfies 
$\int_{E\times [0,1]} m(de) \, l(du) \, \sup_{k\in K}|f(e,k,u)|<\infty$, then
\begin{equation}
\label{cexpl}
\hat \bP\left [\int_{E\times K\times [0,1]} d\Psi^1_t  f \, \Big | \,
\cF_t^{X^1} \right]
=\int_{E\times K \times [0,1]}
 X^1_t(de\times dk) \, l(du) \, f(e,k,u).
\end{equation}
Moreover, if $f(e,k,u)$ satisfies 
$\int_{E\times [0,1]} m(de) \, l(du) \, 
(\exp({\sup_{k\in K}f(e,k,u)})-1)<\infty$, then
\begin{equation}
\label{cexp}
\begin{split}
& \hat \bP\left[
\exp\left(\int_{E\times K\times [0,1]} d\Psi^1_t f \right)
\, \Big | \, \cF_t^{X^1} 
\right] \\
& \quad =
\exp\left(\int_{E\times K}
X^1_t(de\times dk) \,
\log\int_{[0,1]} l(du) \, \exp(f(e,k,u))
\right). \\
\end{split}
\end{equation}

\subsection{Measure--valued, continuum--sites, stepping--stone model}

We emphasize that Assumptions \ref{Huntass} and \ref{finitecap}
are still in force.
Consider $\lambda >1$.  We increase the ``local density'' of particles 
in the above construction by replacing $m$ by $\lambda m$ and 
select the levels $\cU^{\lambda}$ to be i.i.d. uniform on 
$[0,\lambda]$ rather than $[0,1]$.  
Following the construction of $\Psi^1$
and $X^1$ above, define
\begin{equation}
\label{defPsilambda}
\Psi^{\lambda}_t=\sum_{u\in \cU^{\lambda}} 
\delta_{(\hat Z_u(t),\kappa_u(t),u)}
\end{equation}
and
\begin{equation}
X^{\lambda}_t=\frac 1{\lambda} \sum_{u\in \cU^{\lambda}}
\delta_{(\hat Z_u(t),\kappa_u(t))}.
\end{equation}

Note that $\kappa_u$ only depends on locations and types
of particles at levels $v \le u$ and
we can construct $\Psi^{\lambda'}$ from
$\Psi^\lambda$ simultaneously for all
$1 \le \lambda' \le \lambda$
by taking $\Psi^{\lambda'}_t$ 
to be the 
restriction of $\Psi^{\lambda}_t$ to the particles with levels in 
$[0,\lambda']$; 
that is,
\begin{equation}
\Psi^{\lambda'}_t=\sum_{u\in \cU^{\lambda}, \, u\leq\lambda'}
\delta_{(\hat Z_u(t),\kappa_u(t),u)}.
\end{equation}
Consequently, we may carry out the obvious
construction 
to build $\Psi_t^\infty$ and $\cU^\infty$
with levels in $[0,\infty[$.  
The 
initial locations and the levels are such that
\begin{equation}
\Psi^{\infty}_0=\sum_{u\in \cU^{\infty}}
\delta_{(\hat Z_u(0),\kappa_u(0),u)}
\end{equation}
is a Poisson random measure with mean measure given 
by
\begin{equation}
A\times B\times C \mapsto \int_A m(dz) \, \mu (z,B) \, l(C),
\end{equation}
where $l$ is now Lebesgue measure on $[0,\infty[$,
and for $1 \le \lambda < \infty$  
each of the $\Psi_t^\lambda$
can now be defined via (\ref{defPsilambda}) with
$\cU^{\lambda}:=\{u\in \cU^{\infty}:u\leq\lambda \}$.  
The analogue of
(\ref{dualid}) becomes 
\begin{equation}
\label{dualid2}
\begin{split}
&\hat \bP\left[
\sum_{
\begin{array}{c}
u_1,\ldots ,u_n \in \cU^{\lambda}\\
u_i \ne u_j, i \ne j
\end{array}
}
f(\hat Z_{u_1}(t),\kappa_{u_1}(t),\ldots ,\hat Z_{u_n}(t),\kappa_{
u_n}(t)) 
\right] \\
& \quad = 
\lambda^n\int_{E^n}
m^{\otimes n}(d\ee) \,
\HH_t f(\ee), \\
\end{split}
\end{equation}
where $\HH_t$ is defined as in (\ref{coalop}).

Define 
$\cF_t^{\lambda}
=\sigma (X^{\lambda}_s,\Psi^{\infty}_s-\Psi^{\lambda}_s:s\leq t)$, 
and note that, as in 
(\ref{cexpl}) and (\ref{cexp}), we have
\begin{equation}
\label{cexpl2}
\hat \bP\left[\int_{E\times K\times [0,\lambda ]}
d\Psi^{\lambda}_t f
\, \Big | \, \cF_t^{\lambda} \right]
= 
\int_{E\times K \times [0,\lambda]} 
X^{\lambda}_t(de\times dk) \, l(du) \, f(e,k,u)
\end{equation}
and
\begin{equation}
\label{cexp2}
\begin{split}
&\hat \bP\left[
\exp\left(
\int_{E \times K \times [0,\lambda]}
d\Psi^{\lambda}_t f 
\right) 
\, \Big | \, \cF_t^{\lambda}
\right] \\
& \quad =
\exp\left(
\int_{E\times K} X^{\lambda}_t(de\times dk) \, 
\lambda 
\log
\left(
1+\frac 1{\lambda}
\int_{[0, \lambda]}l(du) \, \left(\exp(f(e,k,u))-1 \right)\right)\right).\\
\end{split}
\end{equation}

Suppose that $f(e,k,u)=0$ for $u>\lambda_0$.  Then for $\lambda >
\lambda_0$, 
the random variables on the left of (\ref{cexpl2})  and 
(\ref{cexp2}) do 
not depend on $\lambda$.  Since for fixed $t \ge 0$ the $\sigma$--fields
$\cF^{\lambda}_t$ are decreasing in 
$\lambda$, the left sides of (\ref{cexpl2}) and (\ref{cexp2}) are 
positive, reverse martingales, 
and hence converge $\hat \bP$-a.s. 
as $\lambda \uparrow \infty$. It follows that $X^{\lambda}_t$ converges 
$\hat \bP$-a.s. to a random measure $X^\infty_t$ satisfying
\begin{equation}
\label{cexpl3}
\begin{split}
&\hat \bP\left[\int_{E\times K\times [0,\infty[}d\Psi^{\infty}_t f
\, \Big | \, \cF_t^{X^\infty}\right] \\
& \quad =
\int_{E\times K \times [0,\infty[} 
X_t^\infty(de\times dk) \, l(du) \, f(e,k,u), \\
\end{split}
\end{equation}
and
\begin{equation}
\label{cexp3}
\begin{split}
&\hat \bP\left[
\exp\left(\int_{E\times K\times [0,\infty[}
d\Psi^{\infty}_t f 
\right)
\, \Big | \, \cF_t^{X^\infty} 
\right] \\
& =\exp\left(
\int_{E\times K \times [0,\infty[}
X^\infty_t(de\times dk) \, l(du) \,
\left(\exp f(e,k,u) - 1 \right) 
\right). \\
\end{split}
\end{equation}
In particular, by (\ref{cexp3}),
for each $t\geq 0$, $\Psi^{\infty}_t$ is a doubly 
stochastic Poisson process (that is, a Cox process) with random 
mean measure given by $X_t^\infty \otimes l$.

Dividing both sides of (\ref{dualid2}) by $\lambda^n$ and letting 
$\lambda\rightarrow\infty$, we have
\begin{equation}
\label{dualid3}
\begin{split}
&\hat \bP\left[\int_{(E\times K)^n} 
X^\infty_t(de_1\times dk_1)\cdots X^\infty_t(de_n\times dk_n) \,
f(e_1,k_1,\ldots ,e_n,k_n)\right]\\
& \quad =
\int_{E^n} m^{\otimes n}(d\ee) \, \HH_t f(\ee).
\end{split}
\end{equation}
Note also that
\begin{equation}
X_t^\infty(\cdot\times K)=m.
\end{equation} 
By Remark \ref{samemeasstruct} we can regard the measure  $m(de) \mu(e, dk)$ 
as an element of $\Xi$ (which we will also denote by $\mu$)
and the random measure $X^\infty_t$ as a
$\Xi$--valued random variable.
By Theorem \ref{defX}, $X^\infty_t$ has the same law
as $X_t$ under $\bQ^\mu$ for each $t \ge 0$.  In fact, 
it is not difficult to show that $(X^\infty_t, \, t \ge 0)$
is a Markov process with the same finite--dimensional
distributions as $X$ under $\bQ^\mu$.
We stress, however, that we have only constructed $X^\infty_t$
as an almost sure limit for each fixed $t \ge 0$ rather than as
an almost sure limit in some space of c\`adl\`ag paths.

\section{Dissimilarity for the stepping--stone process $X$}

Suppose in this section that the reference measure $m$ is finite.
Without loss of generality, we can take $m$ to be a probability measure.

\begin{Def}
\label{defdissim}
Consider $\nu \in \Xi$.
For $n=2,3,\ldots$ define the {\em $n^{\mathrm{th}}$--order 
dissimilarity} of $\nu$ to be the quantity
\begin{equation}
D_n(\nu) := \int m^{\otimes n}(d\ee)
\bigotimes_{i=1}^n \nu(e_i) 
\left(\{\kk \in K^n : k_j \ne k_\ell, \, \text{ for all $j \ne \ell$}\}\right).
\end{equation}
Note that $1 \ge D_2(\nu) \ge D_3(\nu) \ge \cdots \ge 0$.
Write $\breve D(\nu):= \sup\{n : D_n(\nu)>0\}$
for the {\em maximal dissimilarity} of $\nu$,
where we set $\sup \emptyset = 1$.
\end{Def}

As we remarked in the Introduction,
it is possible,
by exactly the same argument used in Proposition 5.1 of \cite{Eva97},
to show that if 
$Z$ (and hence also $\hat Z$)
is a symmetric $\alpha$--stable process on the circle $\bT$
with index $1< \alpha \le 2$, then for fixed $t>0$ there
$\bQ^\mu$--a.s. 
exists a random countable subset $\{k_1, k_2, \ldots\}$
of the type--space $K$ such that for Lebesgue almost all $e \in \bT$ 
the probability measure $X_t(e)$ is a point mass at one of the $k_i$.
Indeed, under suitable hypotheses a similar argument should extend to
certain other processes for which points are regular.
It is clear that if $X_t$ also has finite maximal dissimilarity
$\bQ^\mu$--a.s., then the set $\{k_1, k_2, \ldots\}$
is, in fact, finite $\bQ^\mu$--a.s.  Theorem \ref{fiblockfidis}
below provides a sufficient condition for 
the maximal dissimilarity $\breve D(X_t)$
to be finite.  We are able to verify this condition
when $Z$ (and hence also $\hat Z$) is Brownian motion on $\bT$
(see Corollary \ref{coalintofinite}
and the beginning of the proof of Theorem \ref{XlivesinXio}).
We suspect that the condition is also true for L\'evy processes
on $\bT$ for which points are not essentially polar
(see \cite{EvPe98} for an indication that this might be so).

\begin{Def}
\label{defgammaxi}
Observe that if $n' > n$, then 
\begin{equation} 
\left((Z_1^{[n']}, \ldots, Z_n^{[n']}), 
(\gamma_1^{[n']}, \ldots, \gamma_n^{[n']}),   \xi^{[n']}_{|\bN_n}\right)
\end{equation}
has the same distribution as
$(\ZZ^{[n]}, \ggamma^{[n]}, \xi^{[n]})$, where we write
$\xi^{[n']}_{|\bN_n}(t)$ for the restriction of the partition
$\xi^{[n']}(t)$ to $\bN_n$.  Consequently, 
on some probability space 
$(\Omega^{[\infty]}, \cF^{[\infty]}, \bP^{[\infty]})$
there is an
$E^\infty$--valued process $\ZZ$, an
$\bN^\infty$--valued process $\ggamma$,
 and a process $\xi$ taking values
in the space of partitions of $\bN$ such that, in
an obvious notation, 
$((Z_1, \ldots, Z_n), (\gamma_1, \ldots, \gamma_n), \xi_{|\bN_n})$
has the same distribution as
$(\ZZ^{[n]}, \ggamma^{[n]}, \xi^{[n]})$. 
\end{Def}

\begin{Rem}
\label{exchangeabilityfact}
Recall the definition of $\ZZ^{[n]}$, $\CZZ^{[n]}$
and $\xi^{[n]}$ from Remark \ref{ZndefX}.
Let $\CZZ^{[n]\updownarrow}$ and $\xi^{[n]\updownarrow}$
be defined from $\ZZ^{[n]}$ in a similar manner to $\CZZ^{[n]}$
and $\xi^{[n]}$, with the difference that
when two coordinate processes of $\ZZ^{[n]}$ collide, rather than 
the one with the higher index being killed, a colliding particle is
killed at random independently
of the past (with both possibilities equally likely).
It is immediate from the strong Markov property that
$(\ZZ^{[n]}, \xi^{[n]\updownarrow})$ 
has the same distribution as
$(\ZZ^{[n]}, \xi^{[n]})$ for all $n \in \bN$.
Consider $t \ge 0$
and a bijection $\beta:\bN \rightarrow \bN$ 
 and define $\xi^{(\beta)}(t)$,
a random partition of $\bN$, by
$i \sim_{\xi^{(\beta)}(t)} j$ if and only if
$\beta^{-1}(i) \sim_{\xi(t)} \beta^{-1}(j)$. Then
$((Z_{\beta(i)}(0))_{i \in \bN}, \xi^{(\beta)})$
has the same distribution as
$((Z_i(0))_{i \in \bN}, \xi)$.
In particular, for each
$t \ge 0$ the random partition $\xi(t)$ is {\em exchangeable} 
in the sense of Kingman's definition of exchangeable
random partitions (see Section 11 of \cite{Ald85}).
\end{Rem} 

\begin{Theo}
\label{fiblockfidis}
For any $\mu \in \Xi$ and $t \ge 0$, the maximal dissimilarity
$\breve D(X_t)$ under $\bQ^\mu$ is stochastically dominated by
the number of blocks in the partition $\xi(t)$.
In particular, if for some $t > 0$ the partition
$\xi(t)$ has finitely many blocks $\bP^{[\infty]}$--a.s.,
then $\breve D(X_t) < \infty$, $\bQ^\mu$-a.s. 
\end{Theo}

\begin{proof}
Fix a diffuse probability measure $\kappa$ on $K$ (for example,
$\kappa$ could be fair coin--tossing measure).
Given another probability measure $\rho$ on $K$, let
$\rho^\rightarrow$ denote the push--forward of the measure
$\rho \otimes \kappa$ on $K \times K$ by the mapping
$(k,h) \mapsto (k_1,h_1,k_2,h_2,k_3, \ldots)$ from\
$K \times K$ into $K$.  Let $\rho^\leftarrow$ denote
the push--forward of $\rho$ by the mapping
$k \mapsto (k_1, k_3, k_5, \ldots)$ from $K$ into $K$.
Thus the operations $\rho \mapsto \rho^\rightarrow$
and $\rho \mapsto \rho^\leftarrow$ are one--sided inverses
of each other: we have $(\rho^\rightarrow)^\leftarrow = \rho$.
Given $\mu \in \Xi$, define $\mu^\rightarrow, \mu^\leftarrow \in \Xi$
by $\mu^\rightarrow(e) := \mu(e)^\rightarrow$ and
$\mu^\leftarrow(e) := \mu(e)^\leftarrow$.  
Of course,  the operations $\mu \mapsto \mu^\rightarrow$
and $\mu \mapsto \mu^\leftarrow$ are also one--sided inverses
of each other.  Note for any $\mu \in \Xi$
that $\mu^\rightarrow(e)$ is diffuse for
all $e \in E$.

Fix $t \ge 0$.
It is straighforward
to check from the definition in Theorem \ref{defX} that the distribution of
$X_t^\leftarrow$ under $\bQ^{\mu^\rightarrow}$ coincides with
the distribution of $X_t$ under $\bQ^\mu$.  
(This is, of course, what we expect from the stepping--stone
model interpretation: a model that keeps track
of the types for one trait should look the same as a model that keeps
track of the types for two traits if we don't look at 
one of the traits.)
Clearly, 
$D_n(\nu^\leftarrow) \le D_n(\nu)$ for any $n$ and $\nu \in \Xi$,
so $\breve D(\nu^\leftarrow) \le \breve D(\nu)$.  Consequently, $\breve D(X_t)$
under $\bQ^\mu$ is stochastically dominated by $\breve D(X_t)$
under $\bQ^{\mu^\rightarrow}$.

We can use Remark \ref{extnofXdef}
 to compute  multivariate moments of the form
\begin{equation*}
\bQ^\nu
\left[
\left\{D_{n_1}(X_t)\right\}^{a_1} \ldots 
\left\{D_{n_\ell}(X_t)\right\}^{a_\ell}\right], \; 
n_i \in \{2,3,\ldots\}, \, a_i \in \bN, \, 1 \le i \le \ell, \, \ell \in \bN,
\end{equation*}
and discover that they are independent of
$\nu$ within the class of $\nu \in \Xi$
with the property that $\nu(e)$ is diffuse for
all $e \in E$.
Because $0 \le D_k(X_t) \le 1$  for all $k \ge 2$,
the multivariate moment problem for each
of the vectors $(D_{n_1}(X_t), \ldots, D_{n_\ell}(X_t))$
is well--posed and
hence the joint distribution of
$(D_2(X_t), D_3(X_t), \ldots)$ under $\bQ^\nu$
is the same for all such $\nu$.  
Consequently, the distribution of $\breve D(X_t)$
under $\bQ^\nu$ is also the same for all such
$\nu$.  In particular, if
$\lambda \in \Xi$ is defined by $\lambda(e) := \kappa$
for all $e \in E$, then the distributions of $\breve D(X_t)$
under $\bQ^{\mu^\rightarrow}$ and $\bQ^\lambda$ are the same.

Putting the above observations together, we see that 
it suffices to show that $\breve D(X_t)$ under $\bQ^\lambda$
is stochastically dominated by the number of blocks
of $\xi(t)$.

Let $(\tilde L_i)_{i \in \bN}$ be an i.i.d. sequence of
$K$--valued random variables (which we suppose are also
defined on $(\Omega^{[\infty]}, \cF^{[\infty]}, \cP^{[\infty]})$)
that is
independent of $\ZZ$ with $\tilde L_i$
having distribution $\kappa$.  Put $L_i = \tilde L_{\gamma_i(t)}$,
so that $L_i = L_j$ if and only if
$i \sim_{\xi(t)} j$,
$\bP^{[\infty]}$--a.s.  It follows from Remark \ref{exchangeabilityfact}
that the sequence $((Z_i(0), L_i))_{i \in \bN}$ of 
$E \times K$--valued random variables is exchangeable.

Let $\Delta_i$ denote the point mass at $(Z_i(0), L_i)$.
By an extension of the standard reverse martingale
proof of de Finetti's theorem,
as $n \rightarrow \infty$ the sequence
of random probability measures 
$Y_n := n^{-1} \sum_{i=1}^n \Delta_i$
converges $\bP^{[\infty]}$--a.s.
in the weak topology to a random probability measure $Y$
on $E \times K$.
Moreover, if we let $\cI$ denote the permutation invariant
$\sigma$-field corresponding to $((Z_i(0), L_i))_{i \in \bN}$
(that is, $\cI := \bigcap_n \sigma\{Y_n , Y_{n+1}, \ldots\}$), then
we have
\begin{equation}
\label{condexprep}
\bP^{[\infty]}
\left [
\phi
\left (
(Z_1(0), L_1), \ldots, (Z_n(0), L_n)
\right) 
\, | \, \cI 
\right ]
= \int dY^{\otimes n} \phi
\end{equation}
for any bounded Borel function $\phi$ on $(E \times K)^n$.
(See the proof
of Theorem 2.4 of \cite{DoKu96} for the details of this sort
of argument.)

We claim that $Y$ has
the same distribution as $X_t$ under $\bQ^\lambda$
(recall from Remark \ref{samemeasstruct} that
we can identify
$\nu \in \Xi$ with the probability
measure $m(de) \nu(e)(dk)$ on $E \times K$
and
that $\Xi$--valued random variables become
random probability measures on $E \times K$ when thought of in this way).
If $\phi$ is a  bounded Borel function on
$(E \times K)^n$ for some $n \in \bN$, then,
by (\ref{condexprep}),
\begin{equation*}
\begin{split}
&\bP^{[\infty]}
\left[
\int  dY^{\otimes n} \phi 
\right] \\
& \quad = 
\bP^{[\infty]}
\left[
\phi
\left(
(Z_1(0), L_1), \ldots, (Z_n(0), L_n)
\right)
\right] \\
& \quad =
\bP^{[n]}
\left[
\int \bigotimes_{j \in \Gamma^{[n]}(t)} \kappa(dk_j)
\phi
\left(
(Z_1^{[n]}(0), k_{\gamma_1^{[n]}(t)}), \ldots, 
(Z_n^{[n]}(0), k_{\gamma_n^{[n]}(t)})
\right)
\right] \\
& \quad = 
\bP^{[n]}
\left[
\int \bigotimes_{j \in \Gamma^{[n]}(t)} \lambda(Z_j^{[n]}(t))(dk_j)
\phi
\left(
(Z_1^{[n]}(0), k_{\gamma_1^{[n]}(t)}), \ldots, 
(Z_n^{[n]}(0), k_{\gamma_n^{[n]}(t)})
\right)
\right]. \\
\end{split}
\end{equation*}
Comparing this with the equivalent definition 
of $(X, \bQ^\mu)$ in Remark \ref{ZndefX} shows that
$Y$ does indeed have the same distribution as
$X_t$ under $\bQ^\lambda$.

Finally, by (\ref{condexprep}) we have
\begin{equation}
\begin{split}
D_n(Y)
&= 
\bP^{[\infty]}
\left\{
L_i \ne L_j, \, 1 \le i < j \le n 
\, | \, 
\cI
\right\}\\
&\le
\bP^{[\infty]}
\left\{
\exists \ell_1, \ldots, \ell_n :
L_{\ell_i} \ne L_{\ell_j}, \, 1 \le i < j \le n 
\, | \, 
\cI
\right\}\\
&=
\bone
\left\{
\exists \ell_1, \ldots, \ell_n :
L_{\ell_i} \ne L_{\ell_j}, \, 1 \le i < j \le n 
\right\}\\ 
&=
\bone
\left\{
\xi(t) \; \text{has at least $n$ blocks} 
\right\}. \\
\end{split}
\end{equation}
It is thus certainly the case that
$\breve D(X_t)$ under $\bQ^\lambda$ is stochastically dominated by
the number of blocks of $\xi(t)$.
\end{proof}

\begin{Rem}
If $Z$ and $\hat Z$ are both Hunt processes (that is,
if Assumption \ref{Huntass} holds), then the particle representation of
Section  \ref{partconstructX} can be used to give a somewhat
more direct proof of Theorem \ref{fiblockfidis} (note that
Assumption \ref{finitecap} holds because $m$ is a probability measure).
We can sketch the proof as follows.
The set of levels $\cU^\infty$ in the construction of 
Section \ref{partconstructX} is the set of points
of a Poisson random measure on $[0, \infty[$ with Lebesgue
intensity and hence $\cU^\infty$ is discrete.
The dissimilarity $D_n(X_t^\infty)$ is just the conditional 
probability (conditioning on $X_t^\infty$) 
that the particles with the $n$ lowest levels are all of 
different types.  The argument that lead to (\ref{dualid})
establishes that the total number of types 
exhibited by all particles is just
the number of blocks in the corresponding
coalescing Markov labelled partition.
\end{Rem}

\begin{Rem}
In the spirit of the previous remark, it is easy to see that
$\breve D(X_t)$ is almost surely finite for all $t > 0$
when $Z$ (and hence also $\hat Z$) is Brownian motion on the
circle $\bT$ and $m$ is normalised Lebesgue measure.
Once again, we just sketch the argument as a more
quantitative result will follow from
Corollary \ref{coalintofinite} below.
Almost surely, there will exist two particles, say with the
$m^{\mathrm{th}}$ and $n^{\mathrm{th}}$ lowest levels, $m < n$,
such that by time $t$ these two particles have
collided and after the collision the
first particle moved around the circle
clockwise while the second particle moves around
anti--clockwise until they collided again.
The total number of types
exhibited by all particles at time $t$ is then at most $n-1$.
\end{Rem}

\section{Sample path continuity of the stepping--stone process $X$}

Our aim in this Section is to present a sufficient condition
for $X$ to have continuous sample paths (Theorem \ref{ctspathsthm})
and use it to establish that if the migration Markov process
is a L\'evy process 
or a ``nice'' diffusion, then
$X$ has continuous sample paths (Corollary \ref{ctspathsLevy},
Corollary \ref{ctspathsdiffeg} and Remark \ref{appliestodiff}). 
The proof of Theorem \ref{ctspathsthm} is postponed to the
next section.
We emphasize that we are no longer assuming that the reference
measure $m$ is finite.

\begin{Def}
\label{defTe}
For $\ee = (e_1, e_2) \in E^2$ with $e_1 \ne e_2$, let 
$T^\ee := \inf\{t \ge 0 : Z^{e_1}(t) = Z^{e_2}(t)\}$
denote the first time that $Z^{e_1}$ and $Z^{e_2}$ collide.
\end{Def}

\begin{Theo}
\label{ctspathsthm}
Suppose there exists $\varepsilon > 0$ such that for all 
non-negative $\psi \in L^1(m) \cap L^\infty(m)$,
\begin{equation*}
\limsup_{t \downarrow 0} 
t^{-\varepsilon}
\int m^{\otimes 2}(d\ee)
\psi^{\otimes 2} (\ee)
\bP\{T^\ee \le t\}
< \infty.
\end{equation*}
Then $X$ has continuous sample paths $\bQ^\mu$-a.s.
for all $\mu \in \Xi$.
\end{Theo}

\begin{Cor}  
\label{ctspathsLevy}
Suppose that $Z$ is a L\'evy process on
$\bR^d$ or the torus $\bT^d$ for some
$d \in \bN$, and $m$ is Lebesgue measure.
Then $X$ has continuous sample paths $\bQ^\mu$-a.s.
for all $\mu \in \Xi$.
\end{Cor}

\begin{proof}
For $d \ge 2$ we have that $T^\ee = \infty$, $\bP$--a.s.
for $m^{\otimes 2}$--a.e. $\ee$,
and so Theorem \ref{ctspathsthm} certainly gives the result.
In fact, it follows from the remarks at the beginning of
Section 5 in \cite{Eva97} that $X$ evolves deterministically
and continuously in this case.

Now consider the case where $Z$ is $\bR$-valued.  The $\bT$-valued
case is similar and is left to the reader.

Write $(\bar Z, \bar P^z)$ for the L\'evy process that
is the symmetrisation of $Z$.  That is, the distribution
of $\bar Z$ starting at $0$ is the same as that
of $Z' - Z''$, where $Z',Z''$ are two independent copies of
$Z$ both started at $0$.  Put
\begin{equation}
\bar T^0 := \inf\{t \ge 0 : \bar Z(t) = 0\}.
\end{equation}
Then for non-negative $\psi \in L^1(m) \cap L^\infty(m)$,
\begin{equation}
\begin{split}
\int m^{\otimes 2}(d\ee) \psi^{\otimes 2}(\ee) P\{T^\ee \le t\}
& = 
\int m^{\otimes 2}(d\ee) \psi^{\otimes 2}(\ee) 
\bar P^{e_1-e_2}\{\bar T^0 \le t\} \\
& =
\int m(dx) \bar \psi(x) \bar P^x\{\bar T^0 \le t\},\\
\end{split}
\end{equation}
where
$\bar \psi(x) := \int m(dy) \psi(x+y) \psi(y) \in L^1(m) \cap L^\infty(m)$
and we are, of course, using the shift invariance of $m$.

For $\alpha>0$
write $\bar C^\alpha$, $\bar U^\alpha$ and $\bar e^\alpha$
for the $\alpha$-capacity, $\alpha$-resolvent and 
$\alpha$-energy corresponding
to $\bar Z$ (see Sections I.2, II.3 and II.4 of \cite{Ber96} for definitions).  
Because $\bar Z$ is symmetric, these
coincide with the corresponding dual objects.
Write $\Psi$ for the characteristic exponent of $\bar Z$
(see Section I.1 of \cite{Ber96}).  Note that $\Psi$ is real--valued
and non-negative.

Using the convention $\frac{1}{\infty} = 0$, we have 
from Theorems II.7 and II.13 of \cite{Ber96} that
\begin{equation}
\int m(dx) \bar \psi(x) \bar P^x [\exp(-\alpha \bar T^0)]
= \bar{C}^\alpha(\{0\}) \bar U^\alpha \bar \psi(0)
\leq \frac{\bar U^\alpha \bar \psi(0)}{\bar e^\alpha(\{0\})}
\le \frac{\|\bar \psi\|_\infty}{\alpha \bar e^\alpha(\{0\})}.
\end{equation}

By Proposition I.2 of \cite{Ber96},
\begin{equation}
\Psi (z) \leq c z^2, \quad |z|\geq 1,
\end{equation}
for a suitable constant $c$,
and so, for $\alpha \ge 1$,
\begin{equation}
\alpha \bar e^\alpha(\{0\})
=\frac{\alpha}{2\pi}\int_{-\infty}^{\infty}
\frac{1}{\alpha+\Psi (z)} \, dz
\ge \frac{\alpha}{2\pi}\int_{|z| \ge 1}
\frac{1}{\alpha+c z^2} \, dz
\ge c' \alpha^{\frac{1}{2}}
\end{equation}
for a suitable constant $c'$.

Use the inequality 
$\bone_{[0,t]}(x) 
\le e^{-\alpha x} + 1 - e^{-\alpha t} 
\le e^{-\alpha x} + \alpha t$
and take $\alpha = t^{-\frac{2}{3}}$
to get, for $t \le 1$,
\begin{equation}
\int m(dx) \bar \psi(x) \bar P^x\{\bar T^0 \le t\}
\le 
c'' \|\bar \psi\|_\infty \alpha^{-\frac{1}{2}} + \|\bar \psi\|_1 \alpha t
\le 
c^*(\|\bar \psi\|_1 + \|\bar \psi\|_\infty) t^{\frac{1}{3}},
\end{equation}
for suitable constants $c''$ and $c^*$.
Now apply Theorem \ref{ctspathsthm} with $\varepsilon = 1/3$.
\end{proof}

\begin{Cor}
\label{ctspathsdiffeg}
Let $d$ be a metric inducing the topology of the Lusin space $E$.
Write $B(x,r) := \{y \in E : d(x,y) \le r\}$
for the closed ball of radius $r > 0$ centred at $x \in E$
and $S^r := \inf\{t \ge 0 : Z(t) \notin B(Z(0), r)\}$ for the
time taken by $Z$ to first travel distance $r$ from its starting point.
Suppose that there are constants $\alpha, \beta, \gamma >0$
such that
\begin{equation*}
\limsup_{r \downarrow 0} 
r^{-\alpha} \sup_{x \in E}  m(B(x,r)) < \infty
\end{equation*}
and
\begin{equation*}
\limsup_{r \downarrow 0} 
r^{-\gamma} \sup_{x \in E} 
P^x\{S^{r^\beta} \le r\} < \infty.
\end{equation*}
Then $X$ has continuous sample paths $\bQ^\mu$-a.s.
for all $\mu \in \Xi$.
\end{Cor}

\begin{proof}
For non-negative $\psi \in L^1(m) \cap L^\infty(m)$ 
and $\delta > 0$
we have
\begin{equation}
\begin{split}
& \int m^{\otimes 2}(d\ee) \psi^{\otimes 2}(\ee) \bP\{T^\ee \le t\} \\
& \quad =
\int m^{\otimes 2}(d\ee) \psi^{\otimes 2}(\ee) 
     \bone\{d(e_1,e_2) \le \delta\} \bP\{T^\ee \le t\} \\
& \qquad +
\int m^{\otimes 2}(d\ee) \psi^{\otimes 2}(\ee) 
     \bone\{d(e_1,e_2) > \delta\} \bP\{T^\ee \le t\} \\
& \quad \le
\|\psi\|_1 \|\psi\|_\infty \sup_{x \in E} m(B(x,\delta))
+
2 \|\psi\|_1^2 \sup_{x \in E} P^x\{S^{\delta/2} \le t\}. \\
\end{split}
\end{equation}
Take $\delta = t^\beta$ to get that the hypothesis of
Theorem \ref{ctspathsthm} holds with $\varepsilon = (\alpha \beta) \wedge \gamma$.
\end{proof}

\begin{Rem}
\label{appliestodiff}
The above result can be applied to the case where $Z$ is a 
regular diffusion
on $\bR$ in natural scale.  In this case $m$ is the speed
measure and $\hat Z = Z$.  If $m(dx) = a(x) dx$ with
$a$ bounded away from $0$ and $\infty$ and
$d$ is the usual Euclidean metric on $\bR$, then it is not
difficult to see that the conditions of the corollary hold
for $\alpha = 1$, $\beta < 1/2$, and any $\gamma>0$.
We leave the details to the reader.
\end{Rem}

\section{Proof of Theorem \ref{ctspathsthm}}

The proof of Theorem \ref{ctspathsthm} will involve checking
{\em Kolmogorov's criterion} for the sample path continuity of real--valued
processes of the form $(I(X_t;\phi))_{t \ge 0}$ for suitable
$\phi \in C(K)$, where $I(\cdot;\cdot)$
is defined in Definition \ref{defI}.  The proof will be via several lemmas.

\begin{Rem}
\label{swaprem}
In performing the
necessary moment computations for Kolmogorov's
criterion we will need to consider the various orders
in which particles can coalesce in the coalescing system that
define these moments and estimate
the contribution of each possible sequence
of collisions.  We will repeatedly use the fact  
that if, for fixed $i \ne j$,
 we ``swap'' $Z^{e_i}(t)$ and $Z^{e_j}(t)$ immediately after
a stopping time $S$ for
$\ZZ^\ee$ at which $Z^{e_i}(S)=Z^{e_j}(S)$
to form a new process $\tilde \ZZ^\ee$, then $\tilde \ZZ^\ee$
has the same distribution as $\ZZ^\ee$.
More precisely, if we define
\begin{equation}
\tilde Z_i^\ee(t) :=
\begin{cases}
Z_i^\ee(t), & \text{for $t \leq S$},\\
Z_j^\ee(t), & \text{for $t > S$},\\
\end{cases}
\end{equation}
\begin{equation}
\tilde Z_j^\ee(t) :=
\begin{cases}
Z_j^\ee(t), & \text{for $t \leq S$},\\
Z_i^\ee(t), & \text{for $t > S$},\\
\end{cases}
\end{equation}
and
\begin{equation}
\tilde Z_h^\ee := Z_h^\ee, \; h \notin \{i,j\},
\end{equation}
then, by the strong Markov property,
$\tilde \ZZ^\ee$
has the same distribution as $\ZZ^\ee$.
\end{Rem}

\begin{Def}
For $n',n'' \in \bN$
consider $\ee' \in E^{n'}$ and $\ee'' \in E^{n''}$
with $e_1', \ldots, e_{n'}', e_1'', \ldots e_{n''}''$ distinct.
Define a process $\Phi^{\ee'|\ee''}$ 
taking
values in the collection of finite sequences of
two element subsets of 
$\{e_1', \ldots, e_{n'}', e_1'', \ldots e_{n''}''\}$
and 
stopping times $0 = T_0^{\ee'|\ee''} \le 
T_1^{\ee'|\ee''} \le \ldots$ as follows.
For $e_i' \in \{e_1', \ldots, e_{n'}'\}$,
write 
\begin{equation*}
\begin{split}
S_{e_i'}^{\ee'|\ee''} 
& := 
\inf\{t \ge 0 : Z^{e_i'}(t) = Z^{e_j'}(t) \, \text{for some $j \ne i$} \\
& \quad     \text{or} \, Z^{e_i'}(t) = Z^{e_k''}(t)
            \, \text{for some $k$ such that} \, 
            \CZ_k^{\ee''}(s) \ne \dag \, \text{for all} \, s<t\}.\\
\end{split}
\end{equation*}
For $e_i'' \in \{e_1'', \ldots, e_{n''}''\}$ write
\begin{equation*}
\begin{split}
S_{e_i''}^{\ee'|\ee''} 
& := 
\inf\{t \ge 0 : Z^{e_i''}(t) = Z^{e_j'}(t) \, \text{for some $j$}\\
& \quad     \text{or} \, Z^{e_i''}(t) = Z^{e_k''}(t)
            \, \text{for some $k \ne i$ such that} \,
            \,  \CZ_k^{\ee''}(s) \ne \dag \, \text{for all} \, s<t\}. \\
\end{split}
\end{equation*}
Loosely put, 
we are thinking of the particles starting at coordinates of
$\ee'$ as evolving freely without coalescence whereas the particles starting
at coordinates of $\ee''$ are undergoing coalescence
among themselves.
Moreover,
$S_e^{\ee'|\ee''}$ (if it is finite) is the first time
that the particle starting at 
$e$ (where $e$ is either a coordinate of $\ee'$ or $\ee''$)
collides with a ``living''
particle starting at one of the other coordinates. 

Let $R^{\ee'|\ee''} < n'+n''$ denote the 
cardinality of the random set of time points
\begin{equation}
\left\{S_e^{\ee'| \ee''} : 
e \in \{e_1', \ldots, e_{n'}', e_1'', \ldots e_{n''}''\} \;
\text{and} \;
S_e^{\ee' | \ee''} < \infty
\right\}
\end{equation}
 and, if 
$R^{\ee'|\ee''}>0$, write 
$T_1^{\ee'|\ee''} < \ldots < T_{R^{\ee'|\ee''}}^{\ee'|\ee''}$
for an ordered listing of this set.
Put $T_0^{\ee'|\ee''} := 0$ and
$T_\ell^{\ee'|\ee''} := \infty$ for $\ell>R^{\ee'|\ee''}$.

Set $\Phi^{\ee'|\ee''}(T_0^{\ee'|\ee''}) := \emptyset$.  
For $1 \le k \le R^{\ee'|\ee''}$   
write 
\begin{equation}
\{x_k, y_k\} \subseteq \{e_1', \ldots, e_{n'}', e_1'', \ldots, e_{n''}''\},
\end{equation} 
for the $\bP$--a.s.
unique unordered pair such that
$Z^{x_k}(T_k^{\ee'|\ee''}) = Z^{y_k}(T_k^{\ee'|\ee''})$.
By definition, at most one of $x_k$ and $y_k$ belong
to $\{x_1,y_1, \ldots, x_{k-1},y_{k-1}\}$.
Put $\Phi^{\ee'|\ee''}(T_k^{\ee'|\ee''}) 
:= (\{x_1,y_1\}, \ldots, \{x_k, y_k\})$.
Complete the definition of $\Phi^{\ee'|\ee''}$ by setting
$\Phi^{\ee'|\ee''}(t) := \Phi^{\ee'|\ee''}(T_\ell^{\ee'|\ee''})$, where 
$\ell \ge 0$ is such that
$T_\ell^{\ee'|\ee''} \le t < T_{\ell+1}^{\ee'|\ee''}$.  
\end{Def}

\begin{Def}
Given $\ee \in E^n$, $n \in \bN$, with $e_1, \ldots, e_n$ distinct,
define a  process $\Phi^\ee$ taking
values in the collection of finite sequences of
two element subsets of $\{e_1, \ldots, e_n\}$ and 
stopping times $0 = T_0^\ee \le 
T_1^\ee \le \ldots$ by (with a slight abuse) re-using the definitions
of $\Phi^{\ee'|\ee''}$
and $T_k^{\ee'|\ee''}$
with $\ee' = \ee$ and
$\ee''$ the null vector.  That is, all particles
evolve freely with none of them  killed due
to coalescence. Note that if $n=2$, then
$T_1^\ee = T^\ee$, where $T^\ee$ is the first collision time
from
Definition \ref{defTe}.
\end{Def}

\begin{Not}
Given $\ee' \in E^{n'}$ and $\ee'' \in E^{n''}$,
write $\ee':\ee''$ for the {\em concatenation} of these two vectors.
That is, $\ee':\ee' := (e_1', \ldots, e_{n'}', e_1'', \ldots e_{n''}'') \in E^{n'+n''}$.
\end{Not}

\begin{Not}
For $x \in \bR$ write $\lfloor x \rfloor$ for the greatest
integer less than or equal to $x$.
\end{Not}

\begin{Lemma}
\label{onecollisionbd}
For non-negative $\psi \in L^1(m) \cap L^\infty(m)$, $t \ge 0$,  and
$q,n',n'' \in \bN$  
we have
\begin{equation}
\begin{split}
&\int  m^{\otimes n'} \otimes m^{\otimes n''}
(d\ee' \otimes d\ee'')
\psi^{\otimes n'}(\ee') \psi^{\otimes n''}(\ee'')
\bP\{T_q^{\ee'|\ee''} \le t\}\\
& \quad \le 
\int  m^{\otimes n'} \otimes m^{\otimes n''}
(d\ee' \otimes d\ee'')
\psi^{\otimes n'}(\ee') \psi^{\otimes n''}(\ee'')
\bP\{T_q^{\ee':\ee''} \le t\}\\
& \quad \le 
c(n'+n'',q,\psi)
\left(\int
m^{\otimes 2}(d\ee)\psi^{\otimes 2}(\ee)
\bP\{T^\ee \le t\}\right)^{\lfloor\frac{q}{3}\rfloor}\\
\end{split}
\end{equation}
for some constant
$c(n'+n'',q,\psi) $ that depends only on $n'+n''$, 
$q$ and $\psi$.
\end{Lemma}

\begin{proof}  By definition,
$S_e^{\ee':\ee''} \le S_e^{\ee'|\ee''}$ for any
$e \in \{e_1', \ldots, e_{n'}', e_1'', \ldots e_{n''}''\}$,
and so $T_q^{\ee':\ee''} \le T_q^{\ee'|\ee''}$ for all $q$.
It therefore suffices to show that for $q,n \in \bN$
\begin{equation}
\begin{split}
& \int  m^{\otimes n}
(d\ee)
\psi^{\otimes n}(\ee)
\bP\{T_q^\ee \le t\} \\
& \quad \le 
c(n,q,\psi)
\left(\int
m^{\otimes 2}(d\ee)\psi^{\otimes 2}(\ee)
\bP\{T^\ee \le t\}\right)^{\lfloor\frac{q}{3}\rfloor}.\\
\end{split}
\end{equation}

We begin with some notation.
For any sequence of pairs
\begin{equation*}
H=(\{x_1,y_1\},\ldots,\{x_\ell,y_\ell\}), \text {with $x_i,y_i \in E$, 
and $x_i \ne y_i$ for $1 \le i \le |H| := \ell$},
\end{equation*} 
and $t>0$, define an event
\begin{equation}
A_t^H := \{T^{(x_1,y_1)} \leq T^{(x_2,y_2)} \leq \ldots
            \leq T^{(x_\ell,y_\ell)} \le t \}.
\end{equation}
It is easy
to see that $A_t^{H'} \supseteq
A_t^H$ for any subsequence $H'$ of $H$.
Put $\overline H :=\bigcup_{i=1}^\ell\{x_i,y_i \} \subseteq  E$.
For $z \in \overline H$, define
$\iota(z,H)=\{ 1\leq i\leq |H| : z \in \{x_i,y_i\} \}$
to be the set of indices of the pairs in which $z$ appears.

Now fix $G=(\{x_1,y_1\},\ldots,\{x_q,y_q \})$,
with $x_i,y_i \in E$, 
$x_i \ne y_i$, $1 \le i \le q$, and
$\overline G := \bigcup_{i=1}^q\{x_i,y_i \}
\subseteq \{e_1,\ldots,e_n\}$.
We wish to estimate $\bP(A_t^G)$.

Let $(i_1,\ldots,i_h)$ be the subsequence of
$(1,2,\ldots,q)$ obtained by listing the elements of
$\{\max \iota(z,G), \, z \in \overline G\}$
in increasing order. Define a subsequence 
$G_*$ of
$G$ by
\begin{equation}
G_* := 
(\{x_{i_1},y_{i_1}\},\ldots,\{x_{i_h},y_{i_h}\})
=:(\{x_{0,1},y_{0,1}\},\ldots,\{x_{0,|G_*|},y_{0,|G_*|}\} )
\end{equation}
(the reason for the alternative indexing will become
clear as we proceed).
Note that $|\overline G|=|\overline{G_*}| $ because for all $z\in \overline G$,
$z\in \{x_{\max \iota(z,G)}, y_{\max \iota(z,G)}\}\subseteq \overline{G_*}.$
By definition, for 
$1\leq j\leq |G_*|$
the inequalities $\max \iota(x_{i_j}, G_*) \ge j$
and $\max \iota(y_{i_j}, G_*) \ge j$ hold, and at
least one of these inequalities is an equality.
In other words,
\begin{equation}
\min\{\max \iota(x_{0,j},G_*),
\max \iota(y_{0,j},G_*)\}=j, \quad 1\leq j\leq |G_*|.
\end{equation}

Without loss of generality we can assume that
$i_1=\max \iota(x_{i_1},G)\leq \max
\iota(y_{i_1},G)$.
Then
$ x_{i_1}\not\in \{x_r,y_r \}$ for $i_1 <r\leq q $, and, {\em a fortiori},
$ x_{i_1}\not\in \{x_{i_p},y_{i_p} \}$ for $1<p\leq h $.
Hence, $\iota(x_{0,1},G_*)=\iota(x_{i_1},G_*)=\{1\} $ and we are now
in one of the following three cases:

\bigskip
\noindent
{\bf Case I:}  $|\iota(y_{0,1},G_*)|=1$.
Let $G_1$ be the subsequence of $G_*$ 
obtained by deleting $\{x_{0,1},y_{0,1}\}$.
Then 
$ \overline G_1 \cap \{x_{0,1},y_{0,1}\}=\emptyset$,
\begin{equation}
\begin{split}
\bP(A_t^G) 
\leq \bP(A_t^{G_*})
& \leq \bP\left(\{T^{(x_{0,1},y_{0,1})} \le t\} \cap A_t^{G_1}\right) \\
& =P\{T^{(x_{0,1},y_{0,1})} \le t\} \bP(A_t^{G_1}), \\
\end{split}
\end{equation}
and $ |\overline G_1|=|\overline{G_*}|-2=|\overline G|-2 $.

\bigskip
\noindent
{\bf Case II:} $|\iota(y_{0,1},G_*)|=2$.
Write $\iota(y_{0,1},G_*)=\{1,j_2\} $.
Define a subsequence $G_1$ of $G_*$ by deleting
$\{x_{0,1},y_{0,1}\}$ and $\{x_{0,j_2},y_{0,j_2}\}) $ from
$G_* $. 
Then $\overline G_1 \cap \{x_{0,1},y_{0,1}\}=\emptyset$,
\begin{equation}
\begin{split}
\bP(A_t^{G})
\leq \bP(A_t^{G_*})
& \leq \bP\left(\{T^{(x_{0,1},y_{0,1})}\le t\} \cap A_t^{G_1}\right) \\
& = \bP\{T^{(x_{0,1},y_{0,1})} \le t\} \bP(A_t^{G_1}), \\
\end{split}
\end{equation}
and $|\overline G_1| \geq |\overline G_*|-3 = |\overline G|-3$.

\bigskip
\noindent
{\bf Case III:} $|\iota(y_{0,1},G_*)|>2$.
Write $\iota(y_{0,1},G_*)=\{1,j_2,\ldots,j_p\}$
where $1<j_2<\ldots<j_p$
and $y_{0,1}=y_{0,j_2}=\ldots=y_{0,j_p} $.
Then $\max \iota(x_{0,j_2},G_*)=j_2 $ because
$\max \iota(y_{0,j_2},G_*)=j_p>j_2 $.
Let
\begin{equation}
G_{**} := (\{x_{0,1},y_{0,1}' \},\ldots,\{x_{0,|G_*|},y_{0,|G_*|}'\}),
\end{equation}
where
\begin{equation}
y_{0,j}':=
\begin{cases}
y_{0,j}, & \text{if $j\leq j_2$ or $y_{0,j}\neq y_{0,1}$},\\
x_{0,j_2}, & \text{if $j>j_2$ and $y_{0,j}=y_{0,1}$.}\\
\end{cases}
\end{equation}
We then have $\overline G_*=\overline G_{**}$ 
and
by, switching $Z^{x_{0,j_2}} $ and $Z^{y_{0,j_2}} $ at time 
$T^{(x_{0,j_2},y_{0,j_2})}$ in the manner described in Remark \ref{swaprem},
we also have 
$\bP(A_t^{G_*}) = \bP(A_t^{G_{**}})$.
Moreover, $\iota(x_{0,1},G_{**})=1$ because 
$x_{0,1} \notin
\bigcup_{j=j_2}^q \{x_{0,j},y_{0,j}\} =\bigcup_{j=j_2}^q\{x_{0,j},y_{0,j}'\}$.
Now $|\iota(y_{0,1},G_{**})|=2$ and we are in Case II with $G_*$
replaced by $G_{**}$. {F}rom the discussion in Case II we know that
there exists a subsequence $G_1$ of $G_{**}$ such that
$ \overline G_1 \cap \{x_{0,1},y_{0,1}\}=\emptyset,$
\begin{equation}
\begin{split}
\bP(A_t^{G}) 
\leq  \bP(A_t^{G_*})
= \bP(A_t^{G_{**}})
& \leq \bP\left(\{T^{(x_{0,1},y_{0,1})} \le t\} \cap  A_t^{G_1}\right) \\
& = \bP\{T^{(x_{0,1},y_{0,1})} \le t\} \bP(A_t^{G_1}), \\
\end{split}
\end{equation}
and $ |\overline G_1| \geq 
|\overline G_{**}| - 3 = |\overline G_*|-3 = |\overline G|-3 $.

\bigskip

The reduction  procedure that transformed
$G$ into $G_1$ can be repeated at least
$(\lfloor |\overline G|/3 \rfloor - 1)^+$ more times.
That is,
for $0 \le \ell  \le \lfloor |\overline G|/3 \rfloor$, there 
exist sequences
$G_\ell = (\{x_{\ell,1}, y_{\ell,1}\}, \ldots, 
\{x_{\ell,|G_\ell|}, y_{\ell,|G_\ell|}\}$ 
such that $G_0 = G$ and for 
$0 \le \ell \le \lfloor |\overline G|/3 \rfloor - 1$
\begin{equation}
\bP(A_t^{G_\ell})
\leq \bP\left(\{T^{(x_{\ell,1},y_{\ell,1})} \le t\} \cap A_t^{G_{\ell+1}}\right)
= \bP\{T^{(x_{l,1},y_{l,1})} \le t\} \bP(A_t^{G_{l+1}}),
\end{equation}
with
$\overline G_{\ell+1} \cap
\{x_{\ell,1},y_{\ell,1}\}=\emptyset$, 
$\overline G_{l+1} \subseteq \overline G_\ell$,
and $|\overline G_{l+1}| \geq |\overline G_\ell|-3$.

It follows that
\begin{equation}
\label{productbound}
\bP(A_t^{G})
\leq \bP\{T^{(x_{0,1},y_{0,1})} \le t\} \bP(A_t^{G_1})
\leq \ldots
\leq \prod_{i=0}^{\lfloor\frac{|\overline G|}{3}\rfloor-1}
\bP \{T^{(x_{i,1},y_{i,1})} \le t\},
\end{equation}
where the sets
$\{x_{i,1},y_{i,1}\}$, $i=0,1,\ldots,\lfloor |\overline G|/3 \rfloor-1$, are 
pairwise disjoint. 

Write $\bG_t^\ee$ for the set of possible values of
$\Phi^\ee(T_q^\ee)$ on the event $T_q^\ee \le t$.
Note that $|\overline G| > q$ for any $G \in \bG^\ee$
and so the rightmost product in (\ref{productbound})
has at least $\lfloor q/3 \rfloor$ terms.
Therefore, if we let $c(n,q,\psi)$ denote a constant that
only depends on $n,q,\psi$ (but not $t$), 
we have from (\ref{productbound}) that
\begin{equation}
\begin{split}
&\int  m^{\otimes n}(d \ee) 
\psi^{\otimes n}(\ee) \bP\{T_q^\ee \le t\}\\
& \quad =
\int  m^{\otimes n}(d\ee)
\psi^{\otimes n}(\ee)
\sum_{G \in \bG_t^\ee}
\bP\{\Phi^\ee(T_q^\ee)=G, \, T_q^\ee \leq t\}\\
&\quad \leq 
\int m^{\otimes n}(d\ee)
\psi^{\otimes n}(\ee)
\sum_{G \in \bG_t^\ee}
\bP(A_t^{G})\\
& \quad \le
c(n,q,\psi)
\left(\int m^{\otimes 2}(d\ee)
\psi^{\otimes 2}(\ee)
\bP\{T^\ee \le t\}\right)^{\lfloor\frac{q}{3}\rfloor}.\\
\end{split}
\end{equation}
\end{proof}

\begin{Not}
Given $\mu \in \Xi$, define $\mu_t \in \Xi$, $t \ge 0$,
by
\begin{equation}
\begin{split}
I(\mu_t;\phi) 
& := \int m(de) \bP\left[\int \mu(Z^e(t))(dk) \phi(e)(k)\right]\\
& = \int m(de)  \bP\left[\int \mu(e)(dk) \phi(\hat Z^e(t))(k)\right],\\
\end{split}
\end{equation}
so that $I(\mu_t;\phi) = \bQ^\mu[I(X_t;\phi)]$.
\end{Not}

\begin{Lemma}
\label{centredmomentlemma}
Consider $\mu \in \Xi$, $t \ge 0$, $q \in \bN$, and 
$\phi = \psi \otimes \chi$ 
with non-negative $\psi \in L^1(m)\cap L^\infty(m)$
and non-negative $\chi \in C(K)$. There exists
a constant $c(\phi,q)$ that only depends on $\phi$ and
$q$ (and not $\mu$ or $t$) such that
\begin{equation*}
\bQ^\mu
\left[
\left\{
I(X_t;\phi) - I(\mu_t;\phi)
\right\}^{2q}
\right]
\le
c(\phi,q)
\left(\int
m^{\otimes 2}(d\ee)\psi^{\otimes 2}(\ee)
\bP\{T^\ee \le t\}\right)^{\lfloor\frac{q}{3}\rfloor}.
\end{equation*}
\end{Lemma}

\begin{proof}
Given $n', n'' \in \bN$ and vectors 
$\ff' \in E^{n'}$ and $\ff'' \in E^{n''}$
with $f_1', \ldots, f_{n'}'$, $f_1'', \ldots f_{n''}''$ distinct,
write
\begin{equation}
\begin{split}
L^{\ff'|\ff''}(t) 
& := \int 
\bigotimes_{i'=1}^{n'} \mu(Z^{f_{i'}'}(t))(dk_{i'}')
\otimes
\bigotimes_{i''=1}^{n''} \mu(Z^{f_{i''}''}(t))(dk_{i'}'') \\
& \quad \times
\chi^{\otimes n'}\left(k_1', \ldots, k_{n'}'\right)
\chi^{\otimes n''}\left(k_{\gamma_1^{\ff''}(t)}'', 
                   \ldots, k_{\gamma_{n''}^{\ff''}(t)}''\right). \\
\end{split}
\end{equation}
Then, by definition,
\begin{equation*}
\begin{split}
& \bQ^\mu
\left[
\left\{
I(X_t;\phi) - I(\mu_t;\phi)
\right\}^{2q}
\right]
 = 
\sum_{i=0}^{2q} (-1)^i {2q\choose i} 
I^{2q-i}(\mu_t;\phi) \bQ^\mu[I^i(X_t;\phi)] \\
& \quad =
\sum_{i=0}^{2q} (-1)^i {2q\choose i} 
\int m^{\otimes {2q-i}}(d\ee') \otimes
m^{\otimes i}(d\ee'') \psi^{\otimes{2q-i}}(\ee') \psi^{\otimes i}(\ee'')
\bP\left[L^{\ee'|\ee''}(t)\right]. \\
\end{split}
\end{equation*}
Therefore, by Lemma \ref{onecollisionbd}, it suffices to show that
\begin{equation}
\begin{split}
& \sum_{i=0}^{2q} (-1)^i {2q\choose i} 
\int m^{\otimes {2q-i}}(d\ee') \otimes
m^{\otimes i}(d\ee'') \psi^{\otimes{2q-i}}(\ee') \psi^{\otimes i}(\ee'') \\
& \quad \times
\bP\left[L^{\ee'|\ee''}(t) 
\bone\left\{T_p^{\ee'|\ee''} \le t < T_{p+1}^{\ee'|\ee''}\right\}\right] \\
& = 0 \\
\end{split}
\end{equation}
for $0 \le p \le q-1$.

For $\ee' \in E^{2q-i}$
and $\ee'' \in E^i$ with
$e_1', \ldots, e_{2q-i}', e_1'', \ldots, e_i''$ distinct,
write $\bS_{j,h}^{\ee'|\ee''}$, 
$0 \le j \le 2q-i$,
$0 \le h \le i$,
for the collection of subsets of
$\{e_1', \ldots, e_{2q-i}', e_1'', \ldots, e_i''\}$
with exactly $j$ elements from
$\{e_1', \ldots, e_{2q-i}'\}$
and exactly $h$ elements from  $\{e_1'', \ldots, e_i''\}$.
Put
\begin{equation} 
C_{j,h}^i := |\bS_{j,h}^{\ee'|\ee''}| = {2q-i\choose j}{i\choose h}, \quad
0 \le j \le 2q-i, \, 0 \le h \le i.
\end{equation}
It is clear by construction that, recalling the transformation
$H \mapsto \overline H$ from the proof of Lemma \ref{onecollisionbd},
\begin{equation*}
\begin{split}
&\int m^{\otimes {2q-i}}(d\ee') \otimes
m^{\otimes i}(d\ee'') \psi^{\otimes{2q-i}}(\ee') \psi^{\otimes i}(\ee'')
\bP\left[L^{\ee'|\ee''}(t)
\bone\left\{T_p^{\ee'|\ee''} \le t < T_{p+1}^{\ee'|\ee''}\right\}\right] \\
& \quad =
\sum_{j,h}
\int m^{\otimes {2q-i}}(d\ee') \otimes
m^{\otimes i}(d\ee'') \psi^{\otimes{2q-i}}(\ee') \psi^{\otimes i}(\ee'') \\
& \qquad \times
\bP\left[L^{\ee'|\ee''}(t) 
\bone\left\{T_{p}^{\ee'|\ee''} \le t < T_{p+1}^{\ee'|\ee''}\right\}
\sum_{S \in \bS_{j,h}^{\ee'|\ee''}}
\bone\left\{\overline{\Phi^{\ee'|\ee''}}(T_p^{\ee'|\ee''}) = S\right\}\right] \\
& \quad =
\sum_{j,h}
C_{j,h}^i
\int m^{\otimes {2q-i}}(d\ee') \otimes
m^{\otimes i}(d\ee'') \psi^{\otimes{2q-i}}(\ee') \psi^{\otimes i}(\ee'') \\
& \qquad \times
\bP\left[L^{\ee'|\ee''}(t) 
\bone\left\{T_p^{\ee'|\ee''} \le t < T_{p+1}^{\ee'|\ee''},
\;
\overline{\Phi^{\ee'|\ee''}}(T_p^{\ee'|\ee''}) 
= \{e_1', \ldots, e_j', e_1'', \ldots, e_h''\}\right\}\right]. \\
\end{split}
\end{equation*}
Note that 
$|\overline{\Phi^{\ee'|\ee''}}(T_p^{\ee'|\ee''})| \le 2p$,
and so
a necessary condition 
on $j,h$ for a summand in the last term
to be non-zero is that $j+h \le 2p < 2q$.

For fixed $\ee'$ and  $\ee''$ write
$\ee^* := 
(\ee_1', \ldots, \ee_{2q-i}', e_{h+1}'', \ldots, e_i'')$
and
$\ee^{**} := (e_1'', \ldots, e_h'')$.
Observe that
\begin{equation}
\begin{split}
&\left\{T_p^{\ee'|\ee''} \le t < T_{p+1}^{\ee'|\ee''},
\;
\overline{\Phi^{\ee'|\ee''}}(T_p^{\ee'|\ee''}) 
= \{e_1', \ldots, e_j', e_1'', \ldots, e_h''\}\right\} \\
& \quad =
\left\{T_p^{\ee^*|\ee^{**}} \le t < T_{p+1}^{\ee^*|\ee^{**}},
\;
\overline{\Phi^{\ee^*|\ee^{**}}}(T_p^{\ee^*|\ee^{**}}) 
= \{e_1', \ldots, e_j', e_1'', \ldots, e_h''\}\right\}. \\
\end{split}
\end{equation}
Moreover, on this event
the partition $\xi^{\ee^{**}}(t)$ is the restriction of the partition
$\xi^{\ee''}(t)$ to $\bN_h$,
and hence
$L^{\ee'|\ee''}(t) = L^{\ee^* | \ee^{**}}(t)$
on this event.
Therefore, the quantity
\begin{equation}
\begin{split}
&\int m^{\otimes {2q-i}}(d\ee') \otimes
m^{\otimes i}(d\ee'') \psi^{\otimes{2q-i}}(\ee') \psi^{\otimes i}(\ee'') \\
& \quad \times \bP\left[L^{\ee'|\ee''}(t) 
\bone\left\{T_p^{\ee'|\ee''} \le t < T_{p+1}^{\ee'|\ee''},
\;
\overline{\Phi^{\ee'|\ee''}}(T_p^{\ee'|\ee''}) 
= \{e_1', \ldots, e_j', e_1'', \ldots, e_h''\}\right\}\right] \\
\end{split}
\end{equation}
does not vary as $i$ ranges from $h$ to $2q-j$.

The proof is complete once we note that
for fixed $h,j$ with $h < 2q-j$ we have
\begin{equation}
\begin{split}
&\sum_{i=h}^{2q-j}(-1)^i {2q \choose i} C_{j,h}^i
\\
& \quad =\frac{(2q)!(-1)^h}{(2q-j-h)!j!h!}
\sum_{i=h}^{2q-j}
(-1)^{i-h}
\frac{(2q-j-h)!}{(i-h)!(2q-i-j)!}\\
& \quad = \frac{(2q)!(-1)^h}{(2q-j-h)!j!h!}(1-1)^{2q-j-h}\\
& \quad =0.\\
\end{split}
\end{equation}

\end{proof}

\bigskip
\noindent
{\bf Completion of the Proof of Theorem \ref{ctspathsthm}.}
Because $X$ (as a Hunt process) has c\`adl\`ag
paths $\bQ^\mu$-a.s. for all $\mu \in \Xi$,
it suffices to show that $I(X_\cdot;\phi)$ has
continuous sample paths $\bQ^\mu$-a.s. for all $\mu \in \Xi$
and all $\phi$ belonging to some countable subset of $L^1(m,C(K))$
that is separating for $\Xi$.
Moreover, because $I(X_\cdot;\phi)$ already
has c\`adl\`ag paths, verifying Kolmogorov's criterion
establishes that these paths are, in fact, $\bQ^\mu$--a.s. continuous.
That is, verifying Kolmogorov's criterion
does more than just establish the existence of a continuous
version of $X$, it establishes that the version we already have
is continuous.  

Let $\{\hat U^\alpha\}_{\alpha>0}$ denote the resolvent
corresponding to the semigroup $\{\hat P_t\}_{t \ge 0}$.
Suppose that $\cS$ is a countable collection of bounded,
$m$-integrable, continuous, non-negative functions on $E$
with dense linear span in $L^1(m)$ (such a collection
can be seen to exist by combining Lemma A.1 of \cite{Eva97}
with Proposition 3.4.2 of \cite{EtKu86}).
Note that if $\theta \in \cS$, then $\alpha \hat U^\alpha \theta$ converges to
$\theta$ pointwise as $\alpha \rightarrow \infty$.  Also,
$\int m(dx) \alpha \hat U^\alpha \theta(x) = \int m(dx) \theta(x) < \infty$
by the duality hypothesis for the pair $Z, \hat Z$.
By a standard extension of Lebesgue's
dominated convergence theorem (see, for example,
Proposition 18 in Chapter 11 of \cite{Roy68}), if 
$g \in L^\infty(m)$, then
$\lim_{\alpha \rightarrow \infty} \int m(dx) \alpha \hat U^\alpha \theta(x) g(x)
= \int m(dx) \theta(x) g(x)$. 

Write 
$\cD := \{\hat U^\alpha \theta : \theta \in \cS, \, \alpha \, \text{rational}\}
\subseteq L^1(m) \cap L^\infty(m)$.
It follows easily
from what we have just observed that if $\tilde C$ is
a countable dense subset of $\{\chi \in C(K) : \chi \ge 0\}$, 
then the countable
collection of functions of the form $\psi \otimes \chi$,
with $\psi \in \cD$ and $\chi \in \tilde C$, is
separating for $\Xi$.

Fix $\psi \in \cD$ 
(with $\psi = \hat U^\alpha \theta$ for $\theta \in \cS$ and
$\alpha$ rational), $\chi \in \tilde C$, and $q \in \bN$
such that $\lfloor q/3 \rfloor \varepsilon > 1$, where
$\varepsilon > 0$ is as in the statement of the theorem.
In order to show that $I(X_\cdot; \psi \otimes \chi)$
has $\bQ^\mu$-a.s. continuous sample paths for all $\mu$,
it suffices by the Markov property of $X$ and Kolmogorov's
continuity criterion to show for some constants $c$ 
and $\delta$ which
depend only on $\psi, \chi, q$ that
\begin{equation}
\bQ^\mu
\left[
\left\{
I(X_t;\psi \otimes \chi) - I(\mu; \psi \otimes \chi)
\right\}^{2q}
\right]
\le ct^{1+\delta}
\end{equation}
for all $t \ge 0$ and $\mu \in \Xi$.
This, however, follows from Lemma \ref{centredmomentlemma}
and the observation that
\begin{equation}
\begin{split}
& \left|I(\mu_t; \psi \otimes \chi) - I(\mu;\psi \otimes \chi)\right|
\le
\int m(dx) \left|\hat P_t \psi(x) - \psi(x)\right| \\
& \quad = 
\int m(dx) 
\left|
\int_t^\infty \left(e^{-\alpha(s-t)} - e^{-\alpha s}\right)   
\hat P_s \theta(x) \, ds
-
\int_0^t e^{-\alpha s} \hat P_s \theta(x) \, ds
\right| \\
& \quad \le
2\alpha^{-1} (1 - e^{-\alpha t}) \int m(dx) \theta(x)
\le 2 t \int m(dx) \theta(x), \\
\end{split}
\end{equation}
where we have  used the consequence of the
duality hypothesis on $Z$, $\hat Z$ that
$\int m(dx) \alpha \hat P_t \theta(x) = \int m(dx) \theta(x)$.

\section{Coalescing and annihilating circular Brownian motions}
\label{coalannsect}

In this section we develop a duality relationship between 
systems of coalescing Brownian motions on
$\bT$, the circle of circumference $2\pi$, and systems of
annihilating Brownian motions on $\bT$ (Proposition \ref{coalanndual}).
This relation will be used in Section \ref{finmanypure} to investigate
the properties of the stepping--stone model $X$ when the
migration process is Brownian motion on $\bT$.  It will
also be used in Section \ref{coaltree} to study the random tree
associated with infinitely many coalescing Brownian motions on $\bT$.
We mention in passing that coalescing Brownian motion has recently
become a topic of renewed interest (see, for example, \cite{ToWe97}
and \cite{Tsi98}).

For the rest of this paper,
$Z$ (and hence $\hat Z$) will be standard Brownian
motion on  $\bT$,  and
$m$ will be normalised Lebesgue measure on $\bT$.

Given a finite non--empty set $A \subseteq \bT$, enumerate $A$ as
$\{e_1, \ldots, e_n\}$, put $\ee:=(e_1, \ldots, e_n)$, and 
define a process $W^A$, the {\em set--valued coalescing circular
Brownian motion},  taking values in the collection of
non--empty finite subsets of $\bT$ by
\begin{equation}
W^A(t) 
:= \{\CZ_{\gamma_1^\ee(t)}^\ee(t), \ldots, \CZ_{\gamma_n^\ee(t)}^\ee(t)\}
= \{Z_{\gamma_1^\ee(t)}^\ee(t), \ldots, Z_{\gamma_n^\ee(t)}^\ee(t)\},
\quad t \ge 0.
\end{equation}
Equivalently, $W^A(t)$ is the set
of labels of the coalescing Markov labelled partition process 
$\zeta^\ee(t)$.  Of course, different enumerations
of $A$ lead to different processes, but all these processes
will have the same distribution. In words, $W^A$ describes
the evolution of a finite set of indistinguishable
Brownian particles with the feature that particles
evolve independently between collisions but when
two particles collide they coalesce into a single particle.

Write $\cO$ for the  collection of open subsets
of $\bT$ that are either empty or consist of a finite
union of open intervals with distinct end--points. 
Given $B \in \cO$, define 
on some probability space $(\Sigma, \cG, \bQ)$
an $\cO$--valued process
$V^B$, the {\em annihilating circular Brownian motion} as follows.
The end--points of the constituent intervals execute independent
Brownian motions on $\bT$ until they collide, at which
point they annihilate each other.  If the
two colliding end--points are from different intervals, then
those two intervals merge into one interval.
If the two colliding end--points are from the same
interval, then that interval vanishes
(unless the interval was arbitrarily close to $\bT$
just before the collision,
in which case the process takes the value $\bT$).
The process is stopped when it hits
the empty set or $\bT$. 

We have the following duality relation between
$W^A$ and $V^B$.

\begin{Prop}
\label{coalanndual}
For all finite, non--empty subsets $A \subseteq \bT$, all sets $B \in \cO$,
and all $t \ge 0$,
\begin{equation*}
\bP\{W^A(t) \subseteq B\} = \bQ\{A \subseteq V^B(t)\}.
\end{equation*}
\end{Prop}

\begin{proof} 
For $N \in \bN$, let $\bZ_N := \{0,1,\ldots N-1\}$ 
denote the integers modulo $N$.  Let 
$\bZ_N^{\frac{1}{2}} := \{\frac{1}{2}, \frac{3}{2}, \ldots, \frac{2N-1}{2}\}$
denote the half--integers modulo $N$.  A non-empty subset
$D$ of $\bZ_N$ can be (uniquely) decomposed into ``intervals'':
an interval of $D$
 is an equivalence class for the equivalence relation on 
the points of $D$ defined by
$x \sim y$ if and only if
$x = y$, $\{x, x+1, \ldots, y-1, y\} \subseteq D$, or
$\{y, y+1, \ldots, x-1, x\} \subseteq D$ (with all arithmetic modulo $N$).
Any interval other than $\bZ_N$ itself has
an associated pair of (distinct) ``end--points''
in $\bZ_N^{\frac{1}{2}}$: if the interval
is $\{a,a+1, \ldots, b-1, b\}$, then the corresponding
end--points are $a-\frac{1}{2}$ and $b+\frac{1}{2}$
(with all arithmetic modulo $N$). Note that the end--points of
different intervals of $D$
are distinct.

For $C \subseteq \bZ_N$, let $W_N^C$ be a process
on some probability space $(\Omega', \cF', \bP')$ 
taking values in the collection
of non--empty subsets of $\bZ_N$ that is defined in the
same manner as $W^A$, with Brownian motion on $\bT$ replaced
by simple, symmetric (continuous time) random walk on $\bZ_N$
(that is, by the continuous time Markov chain on $\bZ_N$
that only makes jumps from $x$ to $x+1$ or $x$ to $x-1$ at a common 
rate $\lambda>0$ for all $x \in \bZ_N$).
For $D \subseteq \bZ_N$, let $V_N^D$ be a process taking
values in the collection of subsets of $\bZ_N$
that is defined 
on some probability space $(\Sigma', \cG', \bQ')$
in the same manner as $V^B$,
with Brownian motion on $\bT$ replaced by simple, symmetric
(continuous time) random walk on $\bZ_N^{\frac{1}{2}}$ (with the same
jump rate $\lambda$
as in the definition of $W_N^C$).  That is,
end--points of intervals evolve as annihilating random walks
on $\bZ_N^{\frac{1}{2}}$.

The proposition will follow by a straightforward weak limit argument
if we can show the following duality relationship between the
coalescing ``circular'' random walk
$W_N^C$ and the
annihilating ``circular'' random walk
$V_N^D$:
\begin{equation}
\label{approxduality}
\bP'\{W_N^C(t) \subseteq D\} = \bQ'\{C \subseteq V_N^D(t)\}
\end{equation}
for all non--empty subsets of $C \subseteq \bZ_N$,
all subsets of $D \subseteq \bZ_N$, and all $t \ge 0$.

It is simple, but somewhat tedious,
to establish (\ref{approxduality}) by a generator calculation
using the usual generator criterion for duality
(see, for example, Corollary 4.4.13 of \cite{EtKu86}).
However, as Tom Liggett pointed
out to us, there is an easier route.  A little thought
shows that $V_N^D$ is nothing other than the (simple, symmetric)
voter model on $\bZ_N$.  The analogous relationship between
the annihilating random walk and the voter model on $\bZ$
due to \cite{Sch76}
is usually called the {\em border equation} (see Section 2
of \cite{BrGr80} for a discussion and further references).
The relationship (\ref{approxduality}) is then
just the analogue of the usual duality between the voter model
and coalescing random walk on $\bZ$ and it can
be established in a similar manner by Harris's graphical method
(again see Section 2 of \cite{BrGr80} for a discussion and references).
\end{proof}

\begin{Rem}
We have been unable to find an explicit reference
to Proposition \ref{coalanndual} or its analogue for Brownian motion
on $\bR$.  However, if, in the $\bR$--valued
analogue, one considers a limit 
where $A$ approaches
a dense subset of $\bR$, so that the system of coalescing Brownian
motions converges to a coalescing Brownian flow, then the analogous
result for the limiting flow can be found on p18 of \cite{Arr79}.
\end{Rem}

Recall $\ZZ$ and $\ggamma$ from Definition \ref{defgammaxi}.
Define set--valued processes $W^{[n]}$, $n \in \bN$,
and $W$ by 
\begin{equation}
W^{[n]}(t) 
:= \{Z_{\gamma_1(t)}(t), 
     \ldots, Z_{\gamma_n(t)}(t)\} \subseteq \bT, \; t \ge 0,
\end{equation}
and
\begin{equation}
W(t) 
:= \{Z_{\gamma_1(t)}(t), 
     Z_{\gamma_2(t)}(t), \ldots\} \subseteq \bT, \; t \ge 0.
\end{equation}
Thus, $W^{[1]}(t) \subseteq W^{[2]}(t) \subseteq \ldots$
and $\bigcup_{n \in \bN} W^{[n]}(t) = W(t)$.
Recall that $(W^{[n]}(t)_{t \ge 0}$ has the same law as
$(\{Z_{\gamma_1^{[n]}(t)}^{[n]}(t), 
\ldots, Z_{\gamma_n^{[n]}(t)}^{[n]}(t)\})_{t \ge 0}$.
Put $N(t)
:=|W(t)|$, the cardinality of the
random set $W(t)$.
Note that $N(t)$ is also the number of blocks in the partition $\xi(t)$,
which is in turn the cardinality of the random set $\Gamma(t)$.
It is clear that $\bP^{[\infty]}$-a.s. 
$N(t)$ is a non--increasing, right--continuous
function of $t$ and   
if $N(t_0) < \infty$ for some $t_0 \ge 0$, then $N(t) - N(t-)$ 
is either $0$ or $-1$
for all $t>t_0$.
By the following corollary, $N(t) < \infty$, $\bP^{[\infty]}$-a.s., 
for all $t > 0$.

\begin{Cor}
\label{coalintofinite} 
For $t > 0$,
\begin{equation*}
\bP^{[\infty]}\left[\,N(t)\,\right]
=
1 + 2 \sum_{n=1}^\infty 
\exp\left(- \left(\frac{n}{2}\right)^2 t \right)
< \infty
\end{equation*}
and 
\begin{equation*}
\lim_{t \downarrow  0} t^{\frac{1}{2}}
\bP^{[\infty]}\left[\,N(t)\,\right]
=
2 \sqrt{\pi}.
\end{equation*}
\end{Cor}

\begin{proof}
Note that if $B$ is a single open interval
(and hence for all $t \ge 0$ the set
 $V^B(t)$ is either an interval or empty)
and we let $L(t)$ denote the length of $V^B(t)$, then
$L$ is a Brownian motion on $[0,2\pi]$
with $\mathrm{Var} L(t) = 2t$ that is stopped
at the first time it hits $\{0,2\pi\}$.

Now, for $M \in \bN$ and $0 \le i \le M-1$ we have from 
the translation invariance of $Z$ and Proposition
 \ref{coalanndual}
that
\begin{equation}
\begin{split}
& \bP^{[\infty]}\left\{W^{[n]}(t) \cap [2 \pi i / M, 2 \pi (i+1) / M] 
            \ne \emptyset\right\}\\
& \quad = 1 - \bP^{[\infty]}\left\{W^{[n]}(t) \subseteq ]0, 2 \pi (M-1)/M[\right\}\\
& \quad = 1 - \bP^{[\infty]}\left\{W^{[n]}(0) \subseteq V^{]0, 2 \pi (M-1)/M[}(t)\right\},\\
\end{split}
\end{equation}
where we take the 
annihilating process
$V^{]0, 2 \pi (M-1)/M[}$ to be defined on the same probability
space $(\Omega^{[\infty]}, \cF^{[\infty]}, \bP^{[\infty]})$ as
the process $\ZZ$ that was used to construct $W^{[n]}$
and $W$, and we further take the processes $V^{]0, 2 \pi (M-1)/M[}$
and 
$\ZZ$ to be independent.
Thus,
\begin{equation}
\begin{split}
& \bP^{[\infty]}\left\{W(t) \cap [2 \pi i / M, 2 \pi (i+1) / M] \ne \emptyset\right\}\\
& \quad = 1 - 
\bP^{[\infty]}\left\{V^{]0, 2 \pi (M-1)/M[}(t) = \bT\right\}\\
& \quad = 1 - \tilde \bP\left\{\tilde\tau \le 2t, \, 
\tilde B(\tilde\tau) = 2 \pi \; | \; \tilde B(0) = 2 \pi (M-1)/M\right\},\\
\end{split}
\end{equation}
where $\tilde B$ is a standard one--dimensional Brownian motion on 
some probability space $(\tilde\Omega, \tilde\cF, \tilde\bP)$ and 
$\tilde\tau = \inf\{s \ge 0 : \tilde B(s) \in \{0, 2 \pi\}\}$.

By Theorem 4.1.1 of \cite{Kni81} we have
\begin{equation*}
\begin{split}
& \bP^{[\infty]}\left[\,|W(t)|\,\right] \\
& \quad = \lim_{M \rightarrow \infty}
\bP^{[\infty]}\left[\sum_{i=0}^{M-1} 
\bone\left\{W(t) 
\cap [2 \pi i / M, 2 \pi (i+1) / M] \ne \emptyset\right\} \right] \\
& \quad = \lim_{M \rightarrow \infty}
M \left(1 - \tilde\bP\left\{\tilde\tau \le 2t, \, 
\tilde B(\tilde\tau) = 2 \pi 
   \; | \; \tilde B(0) = 2 \pi (M-1)/M\right\} \right) \\
& \quad = 
1 - \lim_{M \rightarrow \infty} M \frac{2}{\pi}
\sum_{n=1}^\infty
\frac{(-1)^n}{n}
\sin\left(n \pi \left(\frac{M-1}{M}\right)\right)
\exp\left(- \left(\frac{n}{2}\right)^2 t \right) \\
& \quad =
1 + 2 \sum_{n=1}^\infty 
\exp\left(- \left(\frac{n}{2}\right)^2 t \right)\\
& \quad = \theta\left(\frac{t}{4\pi}\right) < \infty,\\
\end{split}
\end{equation*}
where
\begin{equation}
\label{deftheta}
\theta(u) := \sum_{n=-\infty}^\infty \exp(- \pi n^2 u)
\end{equation}
is the Jacobi theta function (we refer
the reader to \cite{BiPiYo98} for a survey of
many of the other probabilistic interpretations of
the theta function).  The proof is completed
by recalling that
$\theta$ satisfies the functional equation 
$\theta(u) = u^{-\frac{1}{2}} \theta(u^{-1})$ and
noting that $\lim_{u \rightarrow \infty} \theta(u) = 1$.
\end{proof} 

We conjecture that $t^{\frac{1}{2}} N(t) \rightarrow 2 \sqrt{\pi}$
 as $t \downarrow 0$, $\bP^{[\infty]}$--a.s.  However, we are only able to
prove the following weaker result, which will be used in
Section \ref{coaltree}.  The proof will be given at the end of 
this section after some preliminaries.

\begin{Prop}
\label{asympnoblocks}
With
$\bP^{[\infty]}$--probability one,
\begin{equation*}
0 < \liminf_{t \downarrow 0} t^{\frac {1}{2}} N(t) \le
\limsup_{t\downarrow 0} t^{\frac {1}{2}} N(t) < \infty.
\end{equation*}
\end{Prop}

For $t>0$ the random partition $\xi(t)$ is,
by Remark \ref{exchangeabilityfact} 
and Corollary \ref{coalintofinite},
 exchangeable with a finite number of blocks.
Let $1=x_1^t<x_2^t<\ldots<x_{N(t)}^t$ be the list 
in increasing order of the minimal elements of the
blocks of $\xi(t)$ (that is,
a list in increasing order of the elements of the set $\Gamma(t)$).
  Results of Kingman (see
Section 11 of \cite{Ald85} for a unified account) 
and the fact that $\xi$ evolves by pairwise coalescence of 
blocks give that $\bP^{[\infty]}$--a.s. for all $t>0$
the asymptotic frequencies
\begin{equation}
F_i(t) 
= 
\lim_{n \rightarrow \infty} 
n^{-1}|\{j \in \bN_n : j \sim_{\xi(t)} x_i^t\}|
\end{equation}
exist  for $1 \le i \le N(t)$ and
$F_1(t) + \cdots + F_{N(t)}(t)=1$.

\begin{Lemma}
\label{asympsumsquares}
With
$\bP^{[\infty]}$--probability one,
\begin{equation*}
\lim_{t\downarrow 0}
t^{-\frac {1}{2}} \sum_{i=1}^{N(t)} F_i(t)^2
= 
\frac{2}{\pi^{3/2}}.
\end{equation*}
\end{Lemma}

\begin{proof}
Put $T_{ij} := \inf\{t \ge 0 : Z_i(t) = Z_j(t)\}$ for $i \ne j$.
Observe that
\begin{equation*}
\begin{split}
\bP^{[\infty]}\left[\sum_{i=1}^{N(t)} F_i(t)^2\right]
&=\bP^{[\infty]}\left[
\lim_{n\rightarrow \infty}
\frac{1}{n^2}\sum_{i=1}^n\sum_{k=1}^n 
\bone\left\{j\sim_{\xi(t)} k \right\}\right]\\
&=\bP^{[\infty]}\{1\sim_{\xi(t)} 2\}\\
&=\bP^{[\infty]}\{T_{12}\le t\}.\\
\end{split}
\end{equation*}

{F}rom Theorem 4.1.1 of \cite{Kni81} we have
\begin{equation*}
\begin{split}
&\bP^{[\infty]}\{T_{12} \le t\} \\
& \quad = \frac{1}{2\pi}\int_0^{2\pi}
1 -
\frac{4}{\pi} \sum_{n=1}^\infty
\sin\left(\frac{(2n-1)x}{2}\right)\frac{1}{2n-1}
\exp\left(-\left(\frac{2n-1}{2}\right)^2t\right) \, dx\\
& \quad = \frac{8}{\pi^2}\sum_{n=1}^\infty \frac{1}{(2n-1)^2} 
\left\{1 - \exp\left(-\left(\frac{2n-1}{2}\right)^2t\right)\right\} \\
& \quad = \frac{2}{\pi^2} 
\int_0^t \sum_{n=1}^\infty 
\exp\left(-\left(\frac{2n-1}{2}\right)^2s\right)
\, ds \\
& \quad = \frac{2}{\pi^2}
\int_0^t
\frac{1}{2}
\left\{ \sum_{n=-\infty}^\infty \exp\left(-n^2\frac{s}{4}\right)
- \sum_{n=-\infty}^\infty \exp\left(-n^2 s \right) \right\}
\, ds \\
& \quad = \frac{1}{\pi^2} \int_0^t
\left\{ 
\theta\left(\frac{s}{4\pi}\right) - \theta\left(\frac{s}{\pi}\right)
\right\}
\, ds,\\
\end{split}
\end{equation*}
where
$\theta$ is again
the Jacobi theta function defined in (\ref{deftheta}). 
By the properties of $\theta$ recalled after (\ref{deftheta}),
\begin{equation}
\label{expecgrowth}
\lim_{t \downarrow 0} t^{-\frac{1}{2}} 
\bP^{[\infty]}\left[\sum_{i=1}^{N(t)} F_i(t)^2\right]
= 
\lim_{t \downarrow 0} t^{-\frac{1}{2}}
\bP^{[\infty]}\{T_{12} \le t\}
=
\frac{2}{\pi^{3/2}}.
\end{equation}

Now
\begin{equation*}
\begin{split}
&\bP^{[\infty]}\left[\left(\sum_{i=1}^{N(t)} F_i(t)^2\right)^2\right]\\
&\quad=\bP^{[\infty]}\left[\lim_{n\rightarrow \infty}\frac{1}{n^4} 
\sum_{i_1=1}^n\sum_{i_2=1}^n \sum_{i_3=1}^n\sum_{i_4=1}^n 
\bone\left\{i_1\sim_{\xi(t)}i_2, \, i_3\sim_{\xi(t)}i_4 \right\}\right]\\
&\quad =\bP^{[\infty]}\{1\sim_{\xi(t)} 2, \, 3\sim_{\xi(t)} 4\},\\
\end{split}
\end{equation*}
and so
\begin{equation}
\begin{split}
{\mathrm{Var}}
\left(\sum_{i=1}^{N(t)} F_i(t)^2\right)
&=
\bP^{[\infty]}\{1\sim_{\xi(t)} 2, \, 3\sim_{\xi(t)}4\}
-\bP^{[\infty]}\{T_{12}\le t\}^2 \\
&=
\bP^{[\infty]}\{1\sim_{\xi(t)} 2, \, 3\sim_{\xi(t)}4\}
-\bP^{[\infty]}\{T_{12}\le t, \, T_{23} \le t\}.\\
\end{split}
\end{equation}

Observe that
\begin{equation*}
\begin{split}
&\bP^{[\infty]}
\{T_{12}\le t, \, T_{34}\le t, \, 
   T_{13}>t, \, T_{14}>t, \, T_{23}>t, \, T_{24}>t\}\\
&\quad \le \bP^{[\infty]}\{1\sim_{\xi(t)}2, \, 3\sim_{\xi(t)}4, \, \{\{1,2,3,4\}\} \ne \xi^{[4]}(t)\}\\
&\quad\le  \bP^{[\infty]}\{T_{12}\le t, \, T_{34}\le t\}\\
\end{split}
\end{equation*}
and
\begin{equation*}
\begin{split}
&\bP^{[\infty]}\{T_{12}\le t, \, T_{34}\le t\}
- 
\bP^{[\infty]}
\{T_{12}\le t, \, T_{34}\le t, \, 
   T_{13}>t, \, T_{14}>t, \, T_{23}>t, \, T_{24}>t\}\\
&\quad \le \sum_{i=1,2} \sum_{j=3,4} 
\bP^{[\infty]}\{T_{12}\le t, \,
T_{34}\le t,  \, T_{ij}\le t\}.\\
\end{split}
\end{equation*}
Thus
\begin{equation}
\label{varbound}
\begin{split}
{\mathrm{Var}}
\left(\sum_{i=1}^{N(t)} F_i(t)^2\right)
&\le
\bP^{[\infty]}\{1 \sim_{\xi(t)} 2 \sim_{\xi(t)} 3 \sim_{\xi(t)} 4\} \\
&\quad +
\sum_{i=1,2} \sum_{j=3,4} 
\bP^{[\infty]}\{T_{12}\le t, \,
T_{34}\le t,  \, T_{ij}\le t\}. \\
\end{split}
\end{equation}

Put $D_{ij} := |Z_i(0) - Z_j(0)|$.  We have
\begin{equation}
\label{allequivbd}
\begin{split}
&\bP^{[\infty]}\{1 \sim_{\xi(t)} 2 \sim_{\xi(t)} 3 \sim_{\xi(t)} 4\} \\
& \quad =
\bP^{[\infty]}\{T_{12}\le t, \, T_{13}\wedge T_{23}\le t,
\, T_{14}\wedge T_{24} \wedge T_{34}\le t\} \\
& \quad =
\bP^{[\infty]} 
\biggl(
\{T_{12}\le t, \, T_{13}\wedge T_{23} \le t,
\, T_{14}\wedge T_{24}\wedge T_{34}\le t\} \\
&\qquad \qquad \backslash \;
\{D_{12} \le t^{\frac{2}{5}}, \,
 (D_{13} \wedge D_{23}) \le t^{\frac{2}{5}}, \, 
 (D_{14} \wedge D_{24} \wedge D_{34}) \le  t^{\frac{2}{5}}\}
\biggr) \\
& \qquad +
\bP^{[\infty]}
\{D_{12} \le t^{\frac{2}{5}}, \,
 (D_{13} \wedge D_{23}) \le t^{\frac{2}{5}}, \, 
 (D_{14} \wedge D_{24} \wedge D_{34}) \le  t^{\frac{2}{5}}\} \\
& \quad \le
\sum_{1 \le i < j \le 4} 
\bP^{[\infty]}\{T_{ij} \le t, \, D_{ij} > t^{\frac{2}{5}}\}
+
\bP^{[\infty]}
\left\{\max_{1 \le i < j \le 4} D_{ij} \le 3t^{\frac{2}{5}} \right\},\\
\end{split}
\end{equation}
where we have appealed to the triangle inequality in the last step.
Because $\frac{2}{5} < \frac{1}{2}$,
an application of
the reflection principle and Brownian scaling
certainly gives that the probability
$\bP^{[\infty]}\{T_{ij} \le t, \, D_{ij} > t^{\frac{2}{5}}\}$
is $o(t^\alpha)$ as $t \downarrow 0$ for any $\alpha >0$.
Moreover, by the translation invariance of $m$ 
(the common distribution of the $Z_i(0)$),
the second term in the rightmost member of (\ref{allequivbd}) is at most
\begin{equation*}
\begin{split}
&\bP^{[\infty]}
\{|Z_2(0) - Z_1(0)| \le 3t^{\frac{2}{5}}, \, 
|Z_3(0) - Z_1(0)| \le 3t^{\frac{2}{5}}, \,
|Z_4(0) - Z_1(0)| \le 3t^{\frac{2}{5}}\} \\
&\quad =
\bP^{[\infty]}
\{|Z_2(0)| \le 3t^{\frac{2}{5}}, \, 
|Z_3(0)| \le 3t^{\frac{2}{5}}, \,
|Z_4(0)| \le 3t^{\frac{2}{5}}\} \\
&\quad =
c t^{\frac{6}{5}}, \\
\end{split}
\end{equation*}
for a suitable constant $c$ when $t$ is sufficiently small.
Therefore,
\begin{equation}
\label{varbd1}
\begin{split}
&\bP^{[\infty]}\{1 \sim_{\xi(t)} 2 \sim_{\xi(t)} 3 \sim_{\xi(t)} 4\} \\
& \quad =
\bP^{[\infty]}\{\{T_{12}\le t, \, T_{13}\wedge T_{23}\le t,
\, T_{14}\wedge T_{24} \wedge T_{34}\le t\} \\
& \quad = O(t^{\frac{6}{5}}), \quad \text{as $t \downarrow 0$.} \\
\end{split}
\end{equation}

A similar argument establishes that
\begin{equation}
\label{varbd2}
\bP^{[\infty]}\{T_{12}\le t, \,
T_{34}\le t,  \, T_{ij}\le t\}
= O(t^{\frac{6}{5}}), \quad \text{as $t \downarrow 0$},
\end{equation}
for $i=1,2$ and $j=3,4$.

Substituting (\ref{varbd1}) and (\ref{varbd2}) into (\ref{varbound}) gives
\begin{equation}
{\mathrm{Var}}
\left(\sum_{i=1}^{N(t)} F_i(t)^2\right)
= O(t^{\frac{6}{5}}), \quad \text{as $t \downarrow 0$.}
\end{equation}
This establishes the desired result when combined with the expectation
calculation (\ref{expecgrowth}), Chebyshev's inequality,
a standard Borel--Cantelli argument, and the monotonicity of 
$\sum_{i=1}^{N(t)} F_i(t)^2$.
\end{proof}

We may suppose that on our probability space 
 $(\Omega^{[\infty]}, \cF^{[\infty]}, \bP^{[\infty]})$
there is a sequence $B_1, B_2, \ldots$
of i.i.d. one--dimensional standard Brownian motions 
with  initial distribution  the uniform distribution on $[0,2\pi]$ 
and that $Z_i$ is defined by setting $Z_i(t)$ to be the 
image of $B_i(t)$ under the usual homomorphism from $\bR$ 
onto $\bT$. For $n \in \bN$ and $0 \le j \le 2^n - 1$, let 
$I_1^{n,j} \le I_2^{n,j} \le \ldots$
be a list in increasing order of the set of indices
$\{i \in \bN : B_i(0) \in [{2\pi j}/{2^n},{2\pi(j+1)}/{2^n}[\}$.
Put $B_i^{n,j} := B_{I_i^{n,j}}$ and $Z_i^{n,j} := Z_{I_i^{n,j}}$.
Thus $(B_i^{n,j})_{i \in \bN}$ is an i.i.d. sequence of
standard $\bR$--valued Brownian motions and
$(Z_i^{n,j})_{i \in \bN}$ is an i.i.d. sequence of
standard $\bT$--valued Brownian motions. In each case the
corresponding initial distribution is uniform on
$[{2\pi j}/{2^n},{2\pi(j+1)}/{2^n}[ $.
Moreover, for $n \in \bN$ fixed the sequences 
$(B_i^{n,j})_{i \in \bN}$ are independent as $j$ varies 
and the same is true of the sequences 
$(Z_i^{n,j})_{i \in \bN}$.

Let $\underline {W}$ (resp. $\underline {W}^{n,j}$, $W^{n,j}$)
be the coalescing system defined in terms of 
$(B_i)_{i \in \bN}$ (resp.  $(B_i^{n,j})_{i \in \bN}$,
$(Z_i^{n,j})_{i \in \bN}$) in the same manner that $W$
is defined in terms of $(Z_i)_{i \in \bN}$.

It is clear by construction that
\begin{equation}
\label{comparisonunderl}
N(t) = |W(t)|\leq\sum_{i=0}^{2^n-1}|W^{n,i}(t)|\leq 
\sum_{i=0}^{2^n-1}|\underline{W}^{n,i}(t)|, \quad t>0, \, n \in \bN.
\end{equation}

\begin{Lemma}
\label{finexpunderW} 
The expectation
$\bP^{[\infty]}[ \, |\underline{W}(1)| \, ]$
is finite.
\end{Lemma}

\begin{proof}
 There is an obvious analogue of  
the duality relation Proposition \ref{coalanndual} for systems
of coalescing and annihilating one--dimensional Brownian motions. 
Using this duality and
arguing as in the proof of Corollary \ref{coalintofinite}, 
it is easy to see that,
letting 
$\bar L$ and $\bar U$ be two independent, standard,
real-valued Brownian motions
on some probability space $(\bar\Omega, \bar\cF, \bar\bP)$
with $\bar L(0) = \bar U(0) = 0$, 
\begin{equation*}
\begin{split} 
&\bP^{[\infty]}[|\underline{W}(1)|] \\
&\quad=
\lim_{M \rightarrow \infty} 
\sum_{i=-\infty}^{\infty}
\bP^{[\infty]}\left\{\underline{W}(1) \cap [2 \pi i / M, 2 \pi (i+1) / M] \ne \emptyset\right\}\\
&\quad=
\lim_{M \rightarrow \infty} 
\sum_{i=-\infty}^{\infty}
\bar\bP\biggl\{
\min_{0 \le t \le 1} 
\left( (\bar U(t) + 2 \pi (i+1)/M) - (\bar L(t) + 2 \pi i/M) \right) > 0, \\
& \qquad [\bar L(1) + 2 \pi i/M, \bar U(1) + 2 \pi (i+1)/M] 
\cap [0, 2\pi] \ne \emptyset \biggr\} \\
& \quad \le
\limsup_{M \rightarrow \infty}
c'M \bar\bP\left[
\bone\left\{\min_{0 \le t \le 1} 
\left( \bar U(t)  - \bar L(t) \right) > -2\pi/M\right\} \;
\left(\bar U(1) - \bar L(1) + c''\right)
\right] \\
\end{split}
\end{equation*}
for suitable constants $c'$ and $c''$. 
Noting that $(\bar U - \bar L) / \sqrt{2}$ is a standard Brownian motion,
the result follows from a straightforward calculation with the joint
distribution of the minimum up to time $1$ and value at time $1$ of
such a process (see, for example, Corollary 30 in Section 1.3 
of \cite{Fre83}).
\end{proof}

\bigskip\noindent
{\bf Proof of Proposition \ref{asympnoblocks}.}
By the Cauchy--Schwarz inequality,
\begin{equation} 
1
=
\left(\sum_{i=1}^{N(t)} F_i(t) \right)^2
\le
N(t)  \sum_{i=1}^{N(t)} F_i(t)^2,
\end{equation}
and hence, by Lemma \ref{asympsumsquares},
\begin{equation}
\liminf_{t \downarrow} t^{\frac{1}{2}} N(t) \ge \frac{\pi^\frac{3}{2}}{2},
\quad \bP^{[\infty]}-a.s.
\end{equation}

On the other hand,
for each $n \in \bN$,
$|\underline{W}^{n,i}(2^{-2n})|$, $i=0,\ldots, 2^n-1$, are
i.i.d. random variables which, by Brownian scaling, have 
the same distribution as $|\underline{W}(1)|$.
By (\ref{comparisonunderl}),
\begin{equation}
t^{\frac{1}{2}} N(t)
\leq 
\frac{1}{2^{n-1}} \sum_{i=0}^{2^n-1} 
|\underline{W}^{n,i}(2^{-2n})|
\end{equation}
for $2^{-2n}<t\leq 2^{-2(n-1)}$.
An application of Lemma \ref{finexpunderW}
and the following strong law of large numbers for
triangular arrays completes the proof.

\begin{Lemma}
Consider a triangular array $\{X_{n,i} : 1 \le i \le 2^n, n \in \bN\}$
of identically distributed, real--valued, mean zero, random variables on
some probability space $(\Omega, \cF, \bP)$ such that the collection
$\{X_{n,i} : 1 \le i \le 2^n\}$ is independent for each $n \in \bN$.
Then
\begin{equation*}
\lim_{n \rightarrow \infty} 2^{-n}\left(X_{n,1} + \cdots + X_{n,2^n}\right)
= 0, \; \bP-a.s.
\end{equation*}
\end{Lemma}

\begin{proof}
This sort of result appears
to be known in the theory of complete convergence.
For example, it follows from the much more general
Theorem A in \cite{AsKu80}
by taking $N_n=2^n$ and $\psi(t)=2^t$
in the notation of that result
(see also the Example following that result).  
For the sake of completeness,
we give a short proof that was pointed out to us by Michael Klass.

Let $\{Y_n : n \in \bN\}$ be an independent identically distributed
sequence with the same common distribution as the $X_{n,i}$.
By the strong law of large numbers, for any $\varepsilon > 0$
 the probability that $|Y_1 + \cdots + Y_{2^n}| > \varepsilon 2^n$
infinitely often is $0$.  Therefore, by the triangle inequality,
for any $\varepsilon > 0$ the probability that
$|Y_{2^n+1} + \cdots + Y_{2^{n+1}}| > \varepsilon 2^n$
infinitely often is $0$; and so, by the Borel--Cantelli lemma
for sequences of independent events,
\begin{equation}
\sum_n 
\bP\{|Y_{2^n+1} + \cdots + Y_{2^{n+1}}| > \varepsilon 2^n\} < \infty
\end{equation}
for all $\varepsilon > 0$.
The last sum is also
\begin{equation}
\sum_n
\bP\{|X_{n,1} + \cdots + X_{n,2^n}| > \varepsilon 2^n\},
\end{equation}
and an application of the ``other half'' of the Borel--Cantelli lemma
for possibly dependent events
establishes that for all $\varepsilon > 0$ the probability of
$|X_{n,1} + \cdots + X_{n,2^n}| > \varepsilon 2^n$ infinitely often is $0$,
as required.
\end{proof}

\section{Finitely many pure types for circular Brownian migration}
\label{finmanypure}

Recall that $Z$ and $\hat Z$ are 
standard Brownian motions on the circle $\bT$
and $m$ is normalised Lebesgue measure.  Recall also that
$\cO$ is the collection of open subsets
of $\bT$ that are either empty or the union of a finite number
of disjoint intervals.

\begin{Def}
Let $\Xi^o$ denote the subset of $\Xi$ consisting
 of $\nu$  such that there
exists
a finite set
 $\{k_1^*, \ldots, k_N^*\} \subseteq K$ (depending on $\nu$)
with the property
that for $m$-a.e. $e \in \bT$ we can take $\nu(e) = \delta_{k_i^*}$
for some $i$, and, moreover, we can choose a version of $\nu$
such that the sets 
$\{e \in \bT : \nu(e) = \delta_{k_j^*}\} \in \cO$ for $1 \le j \le N$. 
\end{Def}

\begin{Theo}
\label{XlivesinXio}
For all $\mu \in \Xi$, $\bQ^\mu\{X_t \in \Xi^o \; \text{for all $t>0$}\}=1$.
\end{Theo}

\begin{proof} 
Fix $\mu \in \Xi$ and $t>0$. We will first show that
\begin{equation}
\label{fixedtime}
\bQ^\mu\{X_t \in \Xi^o\}=1.
\end{equation}
 
By the same argument as in Proposition
5.1 of \cite{Eva97}, $\bQ^\mu$-a.s. there is a random countable set of
types $K^*$ such that
$X_t(e) \in \{\delta_k : k \in K^*\}$ for $m$-a.e. $e \in \bT$.
We can also require that $K^*$ has been chosen ``minimally'' so that
$m(\{e \in E : X_t(e) = \delta_k\}) > 0$ for all $k \in K^*$, $\bQ^\mu$-a.s.,
and this requirement specifies $K^*$ uniquely, $\bQ^\mu$-a.s.
For $n \in \bN$ it is clear that on the event where $K^*$ has cardinality
at least $n$ the dissimilarity $D_n(X_t)$ (recall Definition \ref{defdissim})
is strictly
positive $\bQ^\mu$-a.s.
It follows from Theorem \ref{fiblockfidis} and Corollary \ref{coalintofinite}
that $K^*$ is finite 
$\bQ^\mu$-a.s.

In order to show that a representative of 
the equivalence class of $X_t$ in $\Xi$ may be defined so that
$\{e \in \bT : X_t(e) = \delta_k\} \in \cO$
for all $k \in K^*$, it suffices by the device used in the proof
of Theorem \ref{fiblockfidis} to consider the case where
the probability measure $\mu(e)$ is diffuse for all $e \in \bT$
and to show in this case that $\bQ^\mu$-a.s.
for all $k \in K^*$ the
support of the measure $\bone(X_t(e) = \delta_k) \, m(de)$
(which does not depend on the choice of equivalence class representative)
is a connected set.   
For this, it in turn suffices to check that if
$a_1,b_1,c_1,d_1,a_2,b_2,c_2,d_2$ are arranged in anti--clockwise
order around $\bT$, then we have
\begin{equation}
\begin{split}
& \int m^{\otimes 4}(de)
\bone\left\{e_1 \in ]a_1, b_1[, \, e_2 \in ]a_2, b_2[, \,
e_3 \in ]c_1, d_1[, \, e_4 \in ]c_2, d_2[\right\} \\
& \quad \times \int \bigotimes_{i=1}^4 X_t(e_i)(dk_i)
\bone\left\{k_3 \ne k_1 = k_2 \ne k_4\right\}
= 0, \; \text{$\bQ^\mu$-a.s.} \\
\end{split}
\end{equation}
or, equivalently by Remark \ref{extnofXdef},
\begin{equation}
\label{connectedcond}
\begin{split}
& \bone\left\{Z_1^{[4]}(0) \in ]a_1, b_1[, \,
              Z_2^{[4]}(0) \in ]a_2, b_2[, \,
              Z_3^{[4]}(0) \in ]c_1, d_1[,  \,
              Z_4^{[4]}(0) \in ]c_2, d_2[\right\} \\
& 
\quad \times 
\bone\left\{\gamma_3^{[4]}(t) \ne \gamma_1^{[4]}(t) 
            = \gamma_2^{[4]}(t) \ne \gamma_4^{[4]}(t)\right\} = 0,
\; \text{$\bP^{[4]}$-a.s.} \\
\end{split}
\end{equation}

Write, for our fixed $t>0$,
\begin{equation}
T_{ij} = 
\inf \{0 \le s \le t : Z_i^{[4]}(s) = Z_j^{[4]}(s)\}, 
\quad 1 \le i < j \le 4, 
\end{equation}
for the first collision time of
$Z_i^{[4]}$ and  $Z_j^{[4]}$ before time $t$, 
with our standing convention that
$\inf \emptyset = \infty$.
We have 
$\bP^{[4]}\{T_{ij} = T_{k\ell} \ne \infty\} = 0$
for $(i,j) \ne (k,\ell)$.
Suppose that we have a realisation with the properties
\begin{equation}
\label{inbins}
Z_1^{[4]}(0) \in ]a_1, b_1[, Z_2^{[4]}(0) \in ]a_2, b_2[,
Z_3^{[4]}(0) \in ]c_1, d_1[, Z_4^{[4]}(0) \in ]c_2, d_2[,
\end{equation}
\begin{equation}
\label{wronglabels}
\gamma_3^{[4]}(t) \ne \gamma_1^{[4]}(t) 
= \gamma_2^{[4]}(t) \ne \gamma_4^{[4]}(t),
\end{equation}
and
\begin{equation}
\label{distinctexits}
T_{ij} \ne \infty \; \text{implies} \;
T_{ij} \ne T_{k\ell} \; \text{for} \; (i,j) \ne (k,\ell).
\end{equation}

In order that
$\gamma_1^{[4]}(t) = \gamma_2^{[4]}(t)$ holds, 
we must have $T_{12} \ne \infty$.
{F}rom the continuity of the paths
of circular Brownian motion and (\ref{distinctexits}), 
in order that (\ref{inbins}) holds it must then
be the case that
\begin{equation} 
T_{13} \wedge T_{14} 
\wedge T_{23} \wedge T_{24} 
< T_{12} \wedge T_{34}.
\end{equation}
By construction, this would imply that 
$\gamma_3^{[4]}(t)=\gamma_1^{[4]}(t)=\gamma_2^{[4]}(t)$ or 
$\gamma_4^{[4]}(t)=\gamma_1^{[4]}(t)=\gamma_2^{[4]}(t)$, contradicting (\ref{wronglabels}).  Thus
(\ref{connectedcond}) holds and the proof of (\ref{fixedtime})
is complete.

In order to establish the claim of the theorem,
it suffices by (\ref{fixedtime}) and the
Markov property to consider the special case
of $\mu \in \Xi^o$.  Write 
$\{k_1^*, \ldots, k_N^*\} \subseteq  K$ for the corresponding set
of types $k$ such that
$m(\{e \in \bT : \mu(e) = \delta_k\}) > 0$. Fix $1 \le i \le N$.
Let $G \subseteq  K$ be a closed and open
set such that $k_i^* \in G$ and $k_j^* \notin G$ for $j \ne i$
(writing $k_i^* = (h_1,h_2, \ldots)$ one can take
$G= \{(h_1',h_2', \ldots) \in K : h_1' = h_1, \ldots, h_n' = h_n\}$
for some sufficiently large $n$).
It suffices to show for each such $G$ that
if we put $Y_t(e) := X_t(e)(G) \in [0,1]$,
then $\bQ^\mu$--a.s. for all $t \ge 0$ we can choose a representative
of $Y_t \in L^\infty(\bT,m)$ such that 
$Y_t(e) \in \{0,1\}$ for $m$-a.e. $e \in \bT$ and
$\{e \in \bT : Y_t(e) = 1\} \in \cO$.

By the remarks at the end of Section 4 of \cite{Eva97}, we have that
$Y$ is a Feller process with state-space the
subset $L^\infty(\bT,m; [0,1])$
of $L^\infty(\bT,m)$ consisting of $[0,1]$-valued functions
(where $L^\infty(\bT,m; [0,1])$ is equipped with the relative weak$^*$ topology).  
Put $B := \{e \in \bT : \mu(e) = \delta_{k_i^*}\} \in \cO$. 
By the definition of
$X$ in Theorem \ref{defX}
and Proposition \ref{coalanndual}, for $\psi \in L^1(m)$,
\begin{equation}
\begin{split}
& \bQ^\mu\left[\int m^{\otimes n}(d\ee) \psi^{\otimes n}(\ee)
\prod_{i=1}^n Y_t(e_i)\right]\\
& \quad = \int m^{\otimes n}(d\ee) \bP\left[\psi^{\otimes n}(\ee)
\prod_{j \in \Gamma^\ee(t)} \bone_B(Z_j^\ee(t)) \right] \\
& \quad = \int m^{\otimes n}(d\ee) \bP\left[\psi^{\otimes n}(\ee)
\bone\left\{W_t^{\{e_1, \ldots, e_n\}} \subseteq  B\right\}\right] \\
& \quad = \bQ\left[\int m^{\otimes n}(d\ee) \psi^{\otimes n}(\ee)
\bone\left\{\{e_1, \ldots, e_n\} \subseteq  V^B(t) \right\} \right] \\
& \quad = \bQ\left[\int m^{\otimes n}(d\ee) \psi^{\otimes n}(\ee)
\prod_{i=1}^n  \bone_{V^B(t)}(e_i)\right] \\
\end{split}
\end{equation}
(recall that $V^B$ is defined on the probability space
$(\Sigma, \cG, \bQ)$).
Consequently, the $L^\infty(\bT,m; [0,1])$--valued
processes $Y$ and $\bone_{V^B}$ have the same finite--dimensional
distributions.
Clearly, $t \mapsto \bone_{V^B(t)}$ is continuous
(in the weak$^*$ topology).
Therefore,
choosing our representative of $Y_t$ to be $\bone_{V^B(t)}$
for all $t \ge 0$ establishes the desired conclusion.
\end{proof}

\section{The tree associated with coalescing circular Brownian motions}
\label{coaltree}

Recall that $Z$ and $\hat Z$ are standard Brownian motions on
$\bT$ and $m$ is normalised Lebesgue measure.

\begin{Def}
Given
$i,j \in \bN$, let 
$\tau_{ij} := \inf\{t \ge 0 : i \sim_{\xi(t)} j\}$
denote the first time that $i$ and $j$ belong to the same block.
By Remark \ref{exchangeabilityfact},
 the $\tau_{ij}$ are identically distributed. 
Metrise $\bN$ with the (random) metric
$\rho$ given by $\rho(i,j) := \tau_{ij}$.  Observe that
$\rho$ is an {\em ultrametric}; that is, the
{\em strong triangle inequality}
$\rho(i,j) \le \rho(i,k) \vee \rho(k,j)$ holds for all $i,j,k$.
Let $(\bF, \rho)$ denote the completion of $(\bN, \rho)$.
The space $(\bF, \rho)$ is also ultrametric. 
We refer the reader
to Sections 18 and 19 of \cite{Sch84}
for basic facts about ultrametric spaces.
\end{Def}

Some discussion of the space $(\bN, \rho)$ can
be found in Section 4 of \cite{Ald93}.  The analogue
of $(\bF, \rho)$ for another process of coalescing
exchangeable partitions of $\bN$, namely Kingman's
coalescent, is considered in \cite{Eva98} and the
counterpart of Theorem \ref{propsofF} below is obtained.

For background
on Hausdorff and packing dimension see
\cite{Mat95}.
In order to establish some notation, we quickly recall the
definitions of energy and capacity.
Let $(T,\rho)$ be a metric space.  
Write $M_1(T)$ for the collection of (Borel)
probability measures on $T$.
A {\em gauge} is a continuous, non--increasing function
$f:[0,\infty[ \rightarrow [0,\infty]$, such that
$f(r) < \infty$ for $r>0$, $f(0) = \infty$, and
$\lim_{r \rightarrow \infty} f(r) = 0$.
Given
$\mu \in M_1(T)$ and a gauge $f$, the {\em energy
of $\mu$ in the gauge $f$} is the quantity
\begin{equation*}
\cE_f(\mu) := \int \mu(dx) \int \mu(dy) \, f(\rho(x,y)).
\end{equation*}
The {\em capacity of $T$ in the gauge $f$} is the quantity
\begin{equation*}
\Cpty_f(T) := \left(\inf\{ \cE_f(\mu) : \mu \in M_1(T)\}\right)^{-1}
\end{equation*}
(note by our assumptions on $f$ that we need only consider
diffuse  $\mu \in M_1(T)$ in the infimum).

Let $C_{\frac{1}{2}} \subseteq  [0,1]$ denote the middle--$\frac{1}{2}$
Cantor set equipped with the
usual Euclidean metric inherited from $[0,1]$.
One of the assertions of the following result is,
in the terminology of
\cite{PePe95} (see, also, \cite{BePe92, PePeSh96, Per96}),
that $\bF$ is a.s. {\em capacity--equivalent} to 
$C_{\frac{1}{2}}$. Hence, by the results of 
\cite{PePeSh96}, $\bF$ is also a.s. capacity--equivalent to 
the zero set of (one--dimensional) Brownian motion.

\begin{Theo}
\label{propsofF}
With $\bP^{[\infty]}$--probability one, the ultrametric
space $(\bF,\rho)$ is compact with
Hausdorff and packing dimensions both
 equal to $\frac{1}{2}$.
 There exist random variables
$K^*, K^{**}$ such that
$\bP^{[\infty]}$--almost surely
$0 < K^* \le K^{**} < \infty$
and for every gauge $f$
\begin{equation}
\label{capcomparison}
K^* \Cpty_f(C_{\frac{1}{2}}) 
\le \Cpty_f(\bF) 
\le K^{**} \Cpty_f(C_{\frac{1}{2}}).
\end{equation}
\end{Theo}

\begin{proof}
The proof is essentially a reprise of the proof of Theorem 1.1
in \cite{Eva98}, with our Proposition \ref{asympnoblocks} and
Lemma \ref{asympsumsquares} playing the role of
the statements (2.1) and (2.2) in \cite{Eva98}.
\end{proof}

\bigskip
\noindent
{\bf Acknowledgements:} We thank J\"urgen G\"artner, 
Michael Klass, Tom Liggett, Jim Pitman
and Tom Salisbury for helpful discussions.  Some of the research
was conducted while Evans, Fleischmann and Kurtz were visiting the
Pacific Institute for the Mathematical Sciences in Vancouver
and the Mathematical Sciences Research Institute in Berkeley
in Fall 1997.

\bibliographystyle{alpha}

\newpage

\end{document}